\documentclass[a4,12pt]{amsart}

\usepackage{amsmath,amssymb}            
\usepackage{url}

\DeclareMathOperator{\im}{{\sf im}}

\DeclareMathOperator{\Ot}{O}
\DeclareMathOperator{\Sp}{Span}
\DeclareMathOperator{\rk}{rk}

\DeclareMathOperator{\GL}{GL}

\DeclareMathOperator{\id}{id}

\begin{document}
\newcommand{\lab}[1]{\label{#1}}
\renewcommand{\labelenumi}{{\textrm{(\roman{enumi})}}}

\newtheorem{theorem}{Theorem}   
\newtheorem{lemma}[theorem]{Lemma}   
\newtheorem{proposition}[theorem]{Proposition}   
\newtheorem{corollary}[theorem]{Corollary}   
\newtheorem{example}[theorem]{Example}   
\newtheorem{pro}[theorem]{Proposition}
\newtheorem{observation}[theorem]{Observation}
\newtheorem{fact}[theorem]{Fact}
\newtheorem{claim}[theorem]{Claim}

\newcommand{\fdm}{F^{d \times m}}
\newcommand{\fd}{F^{d \times d}}
\newcommand{\ul}[1]{\underline{#1}} 
\newcommand{\pl}{{\sf P}}
\newcommand{\pd}{{\sf D}}
\newcommand{\dv}{{\ul{d}}}
\newcommand{\fdn}{F^{\dv\times\dv}}
\newcommand{\IR}{\mathbb{R}}
\newcommand{\IQ}{\mathbb{Q}}
\newcommand{\IC}{\mathbb{C}}
\newcommand{\IZ}{\mathbb{Z}}
\newcommand{\IN}{\mathbb{N}}

\newcommand{\vep}{\varepsilon} 
\newcommand{\One}{{\bf 1}} 
\newcommand{\Zero}{{\bf 0}} 
  \newcommand{\sk}[2]{\langle #1 | #2 \rangle}
  \newcommand{\mc}[1]{\mathcal{#1}}
  \renewcommand{\phi}{\varphi}

\newcommand{\IL}{\mathbb(L)}
\newcommand{\rr}{{\sf R}}
\newcommand{\vb}{{\sf L}}
\newcommand{\xc}[1]{}

\newcommand{\fforall}{\forall}
\newcommand{\eexists}{\exists}

\newcommand{\bv}{\!\bowtie_{\bar v}\!}
\newcommand{\fn}{(F^{d \times d})^n}

\title[Definable relations]{Definable relations in finite dimensional subspace
lattices with involution. Part II: Quantifier free and homogeneous
descriptions}

\author[Christian Herrmann]{Christian Herrmann}
\address[Christian Herrmann]{Technische Universit\"{a}t Darmstadt FB4\\Schlo{\ss}gartenstr. 7\\64289 Darmstadt\\Germany}
\email{herrmann@mathematik.tu-darmstadt.de}

\author[Martin Ziegler]{Martin Ziegler}
\address[Martin Ziegler]{KAIST School of Computing\\
291 Daehak-ro, Yuseong-gu,\\
Daejeon, South Korea 34141}
\email{ziegler@cs.kaist.ac.kr}


\maketitle

\begin{abstract}
For finite dimensional hermitean inner product spaces $V$,
over  $*$-fields $F$, 
and in the presence of orthogonal bases 
providing form elements in the prime subfield of $F$,
we show that quantifier free definable relations in
 the subspace    lattice $\vb(V)$,
 endowed  with the involution 
induced by orthogonality, 
admit quantifier free descriptions within $F$, also in terms of
 Grassmann-Pl\"ucker coordinates.
In the latter setting, homogeneous
descriptions 
 are obtained if one allows quantification type
$\Sigma_1$.  In absence of involution, these results remain valid.
\end{abstract}

\section{Introduction}

Configurations
in Projective Geometry may be viewed as 
relations on the lattice $\vb(V)$ 
of all subspaces of a vector space $V$, defined by some formula;
e.g.  Desargues's configuration is defined by 
a  conjunction of equations,  dimension conditions
and inequalities to be added.
An orthogonality on $V$ is captured by an involution on $\vb(V)$.

 Back and forth translations of such formulas  
into analytic descriptions, i.e. formulas over the
base field $F$ (with involution) were established
in \cite{hz} for spaces of finite dimension, using the
($*$-)ring of endomorphisms as an intermediate structure.
These translations preserve quantification type $\Sigma_1$
(a block of existential quantifiers followed by a
quantifier free formula); in other words, 
with a $\Sigma_1$-defined relation on 
$\vb(V)$ one associates  the projection of a constructible
affine set over $F$.

In the present note we show that
 quantifier free defined relations on $\vb(V)$
give rise to constructible sets in cartesian as well
as in  Grassmann-Pl\"ucker coordinates. In presence of an
orthogonality we have to require an orthogonal basis
with form elements in the prime subfield.

On the other hand, in the Grassmann-Pl\"ucker setting
with given dimension vector,
we obtain a $\Sigma_1$-formula
 built from homogeneous equations
only, that is, a formula defining a projection of a projective 
constructible subset of the direct product of Grassmannians.
Conversely, any such set can be described within $\vb(V)$
as the projection of a relation defined by
the conjunction of equations $t_i=0$ and $s_j=1$.

In \cite{hz} it was  shown that  preservation
of  quantifier freeness is not possible, in general, for 
translation from $F$ to $\vb(V)$.
 A sufficient condition
for preservation 
 is that  coordinate systems
(or, in the absence of
orthogonality, systems of elements in general position, cf. \cite{gelfand})
are implicitly given by the defining formula:
in analogy to Cayley factorization, cf. \cite{Whiteley,White}.

All translations in this note are effective,  a detailed discussion
of complexity 
shall be postponed to subsequent work. The particular case, where
$F$ is a $*$-subfield of $\IC$ and $V$ admits an orthonormal basis,
has been  studied in \cite{jacm} in the context
of complexity of real computation.

Following the referee's good advice,  we have included
two sections on examples.

\xc{
\section{Examples}
As mentioned above, definable relations in lattices $\vb(V)$
of subspaces  are also  known as ``configurations''.
Of course, definable relations are invariant under automorphisms so that
``configurations'' can be seen both as 
collections of isomorphism types and as instances thereof.

A large supply of examples of quantifier free definable
subsets of $\vb(V)^n$ is given as follows:
 Fix a  partial lattice $P$
of finite  cardinality $n$ and let $M$
consist  of all $(\gamma(p)\mid p \in P)$ where
$\gamma:P \to \vb(V)$ is a homomorphism.  This may be modified by
prescribing $\dim \gamma(p)$ for some $p$
or requiring, for certain $p,q$, the images  
 $\gamma(p),\gamma(q)$ to be related 
in terms of equality, incidence, and their negations;
also in terms of orthogonality if $\vb(V)$ is endowed with an 
involution.

If $P$ is just a partially ordered set then
 configurations $(\gamma(p)\mid p \in P)$ in $\vb(V)^n$
correspond to  linear representations of $P$ in $V$.
According to Corollary~\ref{homcor}, below,
the set of these configurations can 
be defined by homogeneous quantifier free formulas
in terms of  Grassmann-Pl\"ucker coordinates.
A related construction,
for the quiver associated with $P$,
 are quiver Grassmannians (cf. \cite{ring}).

A well known example,
given by a partial lattice,  is the Pasch-Veblen-Young configuration.
Here, $P$ consists of 
 $6$ ``points''  $p_i, p_{i,i+1}$ 
  and $4$  ``lines'' $l$ and  $l_{i,i+1}$, indices modulo $3$.
Points $p_{0,1},p_{1,2},p_{2,3}$ are incident with $l$,
points $p_i$, $p_{i+1}$, $p_{i,i+1}$ are incident with $l_{i,i+1}$ but
$p_{i+2}$, $p_{i+1,i+2}$, $p_{i+2,i+3}$ are not
incident with $l_{i,i+1}$. 
The dimension requirements are $\dim \gamma(p)=1$ for points,
$\dim \gamma(l)=2$ for lines and the maps $\gamma$ have to be injective.
This configuration illustrates a basic axiom of
Projective Geometry and provides coordinate systems for
projective planes. 

Dropping the injectivity and  dimension requirements,
one has to add one more element $b$ to 
the partial lattice and to require 
$l_{i,i+1}\cap p_{i+2}=b$ and
$x+y=z$, $x \cap y=b$
for any points $x\neq y$ incident with line $z$.
 The associated  configurations
in lattices of subspaces are  singletons and 
normalized frames of order $3$, introduced
by von Neumann as coordinate systems for  continuous 
geometries. Singletons can be excluded by 
requiring $\gamma(x)\neq\gamma(y)$ for some $x \neq y$ 
in the partial lattice.

Another simple example is related to the Fano axiom.
The
partial lattice associated with the Pasch configuration 
is extended by one more point $p$
and $3$ more lines $l_i$ such that
$p_i$, $p$, $p_{i+1,i+2}$ are incident with $l_i$, but no other points
are incident with $l_i$;
and such that $p$ is not incident with any of the lines of
the Pasch configuration. Thus, these points and lines form a projective
plane of order $2$. 

Dropping dimension and injectivity requirements,
one adds elements $b,t$ and the lattice relations $x+y=t$, $x \cap y=b$
for each point $x$ non-incident with line $y$ and
$x+y=z$ for any points $x\neq y$ incident with line $z$. 
The associated configurations in a subspace lattice 
are singletons and sublattices isomorphic to the subspace lattice
of a projective plane of order $2$.

On the other hand, replacing the incidence $p_{20} \leq l$ by the non-incidence
condition $\gamma(p_{20}) \not\leq \gamma(l)$ one obtains
configurations failing the Fano axiom. These  give rise to
a subset $M$ of $\vb(V)^{14}$ definable by quantifier free $\phi(\bar x)$.
According to Corollary \ref{p7c}, below,
 $\phi$ translates  into a quantifier free formula
$\tau(\phi)(\bar X)$ defining   the set
  $\theta_{\bar v}(M)$ of $n$-tuples $\bar A$ of matrices
$A_k=(a_{kij})_{ij} \in F^{d \times d}$
such that $\bar u \in M$ where 
$u_k={\rm Span}\{ \sum_{i=1}^d a_{kij} v_i\mid j=1,\ldots, d   \}$;
here $\bar v$ is a fixed basis of $V$ but 
$\theta_{\bar v}(M)$ does not depend on the choice of
$\bar v$.  Now, validity of the Fano axiom means
$M=\emptyset$, that is validity in $F$  of 
the universal sentence
$\forall \bar X\,\bigl(\tau(\phi)(\bar X) \Rightarrow x_{111}\neq x_{111}\bigr)$.
Of course, there are easier ways to do this.
This kind of reasoning together with Theorem \ref{p7}
and the fact that universal sentences 
can be captured by finite systems of forbidden finite partial substructures
yields the following.

\begin{corollary}
Universal sentences in synthetic Projective Geometry
have universal counterparts in coordinates $($dimension being fixed$)$.
\end{corollary}
 }

\section{Preliminaries}\lab{pre}

Statements presented as ``Fact'' 
are well known or obvious; proofs will be omitted or sketched. 
In the sequel, let $F$ 
be a field with involution $r \mapsto r^*$ (and prime subfield $F_0$)
and $V$ a (right) $F$-vector space of (fixed) $\dim V=d <\infty$ 
turned into an inner product space by a
non-alternate  non-degenerate  $*$-hermitean form $\sk{.}{.}_V$
(we will speak just of a \emph{form} and write
$\sk{.}{.}$ if there is no confusion),
that is: additive in both arguments and
\[ \sk{vr}{ws}= r^*\sk{v}{w}s,\;\;
\sk{w}{v}= \sk{v}{w}^*\]
as well as  $\sk{v}{v}\neq 0$ for some $v$, and
$\sk{w}{v}=0$ for all $w \in V$ only if $v=0$  cf. \cite[Chapter I]{gro}.
We write $|v|=\sk{v}{v}$.

A basis  $\bar v=(v_1, \ldots ,v_d)$ of $V$
is  \emph{orthogonal} if $\sk{v_i}{v_j}=0$
for  $i \neq j$;  we will speak of a $\perp$-\emph{basis}.
Recall that such always exist \cite[II \S2 Corollary 1]{gro}: any 
$v_1\neq 0$ can be
completed to a $\perp$-basis.
Given a $\perp$-basis $\bar v$ of $V$, forms are in 1-1-correspondence
with \emph{form constants}  $\bar \alpha=(\alpha_1, \ldots ,\alpha_d) \in F^d$ 
such that $\alpha_i=\alpha_i^*\neq 0$ for all $i$:
namely,   $\alpha_i=|v_i|$ and
 \[\sk{\sum_i v_ir_i}{\sum_j v_js_j}=\sum_k r^*_k \alpha_ks_k.\]
A \emph{scaled isometry} (with factor $r \neq 0$ in $F$)
between inner product spaces $V$ and $W$ over $F$ is 
a linear isomorphism $\omega:V \to W$
such that $\sk{\omega x}{\omega y}_W=r \sk{x}{y}_V$
for all $x,y \in V$. 
Given $\perp$-bases $\bar v$  and
$\bar w$ of  $V$ and $W$ the
linear  isomorphism 
matching $\bar v$ and $\bar w$ is a scaled  isometry
with factor $r$ iff $(|w_1|,  \ldots ,|w_d|)= r(|v_1|,\ldots ,|v_d|)$.
In particular,  $|\bar v|:=(1,|v_1|^{-1}|v_2|, \ldots ,|v_1|^{-1}|v_n|)$
determines the isometry type of $(V,\bar v)$
up to scaling.
Call $\bar \alpha=(\alpha_1, \ldots, \alpha_d)$ \emph{admissible}
(w.r.t. $V$ )
if $\alpha_1=1$ and  $\alpha_i=\alpha_i^*\neq 0$ for all $i$
and 
 if there is
a $\perp$-basis $\bar v$ of $V$ such that
$\bar \alpha=|\bar v|$.

$F^d_{\bar \alpha}$ will denote the space $F^d$ of column vectors
$\ul{u}$
with canonical basis $\bar v$ such that $|v_i|=\alpha_i$.
In this case, $\bar \alpha$ is \emph{admissible}
 if $\alpha_1=1$.

\medskip
The linear subspaces of $V$ form a lattice $\vb(V)$ with bounds
$0,V$, joins given as $U_1+U_2$, meets as $U_1 \cap U_2$,
 and involution 
\[U \mapsto U^ \perp=\{x \in V\mid \forall u\in U.\, \sk{x}{u}=0\}.\]
If $d \geq 3$, then any involution
(that is, an order reversing map  of order $2$) on the lattice
$\vb(V)$ is induced by some kind of  inner product.
Observe that \[U_1 \cap U_2 = (U_1^\perp +U_2^\perp)^\perp\mbox{ and } V^\perp=0.\] 
$\vb(V)$ remains unchanged under scaling the form on $V$.
Moreover, any scaled isometry $\omega:V\to W$ induces
an isomorphism  $\vb(V) \to \vb(W)$: $U \mapsto \omega(U)$.
In particular, given a $\perp$-basis $\bar v$ of $V$ one has an
isomorphism 
 $\Omega_{\bar v}:\vb(V)\to \vb(F^d_{\bar  \alpha})$ where $\bar \alpha=|\bar v|$. If there is no confusion we write $\vb$ in place of $\vb(V)$.

\medskip
 For $m \leq 2d$, we consider  $F^{d \times m}$ the space of
$d\times m$-matrices  $A=(a_{ij})_{ij} $  over $F$
 with columns $\ul{a}_j$.
Let $\rk(A)$ denote the rank of $A$ and
$\Sp(A)$ the $F$-linear subspace of $F^d$ spanned by the columns of $A$;
recall that $B\in \Sp(A)$ iff $B=AC$ for some $C\in F^{d \times d}$.

Write $M_m=\{1,\ldots ,m\}$;
let   $d^\#$ consist of all strictly monotone
maps $f:M_k\to M_d$ where  $k=:|f|\leq d$.
The matrix $A$  is in (reduced) column echelon  \emph{normal form}, shortly NF,
with \emph{positions} $f\in d^\#$ of \emph{pivots} $a_{f(j)j}$ 
(shortly $f$-NF)
if for all $j\leq |f|$,  $h\leq m$, and $i \leq d$
\[\begin{array}{rl}  
 a_{f(j)j}= 1 &\\
 a_{ij}=0 &\mbox{ if } i<f(j)\\
a_{i h}=0 &\mbox{ if } i< f(j+1) \mbox{ and } h>j<|f|\\
a_{f(j) h}=0 &\mbox{ if } h \neq j.
\end{array}\]
$A$ is in \emph{weak normal form} (shortly wNF resp. $f$-wNF) 
if $r A$ is in NF resp. $f$-NF for some $r \neq 0$.
  
$A^*=(a_{ji}^*)_{ij}$ is the conjugate
(w.r.t. the involution on $F$) transpose of $A=(a_{ij})_{ij}$,
$\bar \alpha=(\alpha_1, \ldots ,\alpha_d)\in F^d$ 
$\alpha_i=\alpha_i^*\neq 0$,
$D_{\bar  \alpha}$  the diagonal matrix with diagonal entries
$\alpha_i$. Let $I_k$ denote the unit matrix in $F^{k\times k}$.  
\begin{fact}\lab{nf}
\begin{enumerate}\item
For each $A\in F^{d \times m}$ there is unique
$A^\#\in F^{d \times d}$ in NF such that $\Sp(A)=\Sp(A^\#)$.
Moreover, $|f|=\rk(A)$ iff $A^\#$ in $f$-NF.
For suitable $r\neq 0$, $rA^\#$ can be
obtained from $A$ by  Gaussian   column transformations
without inversion of scalars.
\item For idempotent $A \in F^{d \times d}$
one has $\Sp(A)^\perp=\Sp(I-A^{\dagger_{\bar  \alpha}})$ in $F^d_{\bar  \alpha}$ 
where $A^{\dagger_{\bar  \alpha}}=D_{\bar  \alpha}^{-1}A^*D_{\bar  \alpha}$.
\item For any $f\in d^\#$
 and  $A \in \fd$ in $f$-wNF with 
$A=r A^\#$ one has
$\Sp(A)^\perp=\Sp(P^{-1}(r I- (PA)^{\dagger_{\beta}}))$
where $P$ is a permutation matrix such that
\[\frac{1}{r}PA=\left( \begin{array}{cc}I_{|f|}&0\\\frac{1}{r} K&0
 \end{array} \right) \]
for some $K$ 
and  where  $\bar \beta^t=P\bar \alpha^t$.
\end{enumerate}
\end{fact}
\begin{proof} (i) To avoid inversion,
multiply any column, to be changed, first by a suitable scalar.
Once  echelon form is obtained with each pivot the only non-zero
entry in its row, multiply each pivot with the
product of the others.   For uniqueness see e.g. \cite{yust}.

(ii) This is well known in the context of the  $*$-regular 
ring defined by the form on $F^d_{\bar  \alpha}$. For convenience,
we recall the proof.
 Observe that the linear  map $\phi$ defined 
by $A$ w.r.t. the canonical basis $\bar v$ has adjoint $\phi^\dagger$
in $F^d_{\bar  \alpha}$ defined by $A^{\dagger_{\bar  \alpha}}$.  
Indeed, we have
$\sk{\ul{a}_j}{v_k}=a_{kj}^*\alpha_k$ and
$\sk{v_j}{A^{\dagger_{\bar  \alpha}} v_k}=\sk{v_j}{\alpha_j^{-1}a_{kj}^*\alpha_k v_j}
=\sk{v_j}{v_j}\alpha_j^{-1}a_{kj}^*\alpha_k=a_{kj}^*\alpha_k$.
Now, since $\phi$ is idempotent, so is $\phi^\dagger$
whence $\sk{\phi x}{y-\phi^\dagger y}
=\sk{x}{\phi^\dagger(y-\phi^\dagger y)}=\sk{x}{0}=0$
for all $x,y \in F^d$;
that is $\im \phi =\Sp(A)$ orthogonal to $\im(\id -\phi^\dagger)=
\Sp(I-A^{\dagger_{\bar  \alpha}})$. Since $A^{\dagger_{\bar  \alpha}}$ is idempotent, one has
$\rk(I-A^{\dagger_{\bar  \alpha}})=d-\rk(A^{\dagger_{\bar  \alpha}})=\rk(A)$ whence 
$\Sp(A)^\perp=\Sp(I-A^{\dagger_{\bar  \alpha}})$.

 (iii)
Evidently, 
 $\frac{1}{r}PA$ is idempotent.
By  (ii) applied to the space $F^d_{\bar \beta}$
(with orthogonality $\perp_{\bar \beta}$)
one obtains
$\Sp(PA)^{\perp_{\bar \beta}}=\Sp(\frac{1}{r}PA)^{\perp_{\bar \beta}}= \Sp(I-(\frac{1}{r}PA)^{\dagger_{\bar \beta}})= 
\Sp(r I - (PA)^{\dagger_{\bar \beta}})$.
(iii) follows
observing that $v \mapsto Pv$
defines a linear isometry of $F^d_{\bar \alpha}$
onto $F^d_{\bar \beta}$. 
\end{proof}

\medskip
We consider  first order languages with countably many variables:
$\Lambda_L$ in the signature of $+,0,^\perp$ of bounded
lattices with involution (defining meet by $s \cap t=(s^\perp +t^\perp)^\perp$
and $1=0^\perp$),
$\Lambda_F$ 
in the signature of $*$-rings, with constants $0,1$,
and $\Lambda_F^+$
with  additional constants $c_1, \ldots ,c_d$
to be interpreted as $\alpha_1, \ldots, \alpha_d$
given any $\bar  \alpha$ admissible w.r.t. $V$,
that is, considering $F$ with additional constants
$\alpha_i$.  
A $\Lambda_F^+$-term $t(x_1,\ldots ,x_n)$  is \emph{basic} if it is
a multivariate $*$-polynomial in variables $x_1, \ldots ,x_n$
with  coefficients from $R=\mathbb{Z}[c_1,\ldots ,c_d]$,
that is
an $R$-linear combination of terms $\prod_i x_i^{k_i}(x_i^*)^{l_i}$.
A $\Lambda_F^+$-formula is \emph{basic} 
if it is  a conjunction
of formulas  $p_i=0$ and $q_j\neq 0$
where the $p_i$ and $q_j$ are  basic terms.

An $n$-ary relation $R$  on  $\vb$ resp. $F$ 
is \emph{definable} if there is
a formula $\phi(\bar x)$ (in the relevant language)
such that $R$ consists of all $\bar a$
such that $\phi(\bar a)$ holds in $\vb$ resp. $F$.
Here,
finite strings of variables or elements are written e.g.
as $\bar x$ and $\bar a$, the length being given by  context.
We also use matrices $X=(x_{ij})_{ij}$ of variables in an
obvious way.
$\Sigma_k$ ($\Pi_k$) consists of the prenex formulas with
at most $k$ blocks, each consisting  of quantifiers of the same type, 
the first block being of type $\exists$ ($\forall$).

\begin{fact}\lab{def}
Given a $*$-field $F$  and
values $\alpha_i$ of constants $c_i$.
\begin{enumerate}
\item 
If the $\alpha_i$ are $\Lambda_F$-definable within $F$
then any relation on $F$ which is definable
within $\Lambda_F^+$ is also definable within $\Lambda_F$.
\item Any $\Lambda_F^+$-term is equivalent to a basic term
(uniformly for all $F$ and $\alpha_i$).
\item Any relation defined in $F$
by a quantifier free  $\Lambda_F^+$-formula is
the disjoint union of relations defined by basic
$\Lambda_F^+$-formulas. Moreover, if all $\alpha_i \in F_0$,
then the latter can be chosen in $\Lambda_F$
(that is, with integer coefficients). 
\item If $^*$ is the identity involution of $F$,
any multivariate $*$-polynomial can be considered
a multivariate polynomial in the same coefficients and variables.
\end{enumerate}
\end{fact}

\begin{proof}
(i)
Given a $\Lambda_f^+$-formula $\phi$, 
if the $\alpha_i$ are defined by the
$\alpha_i^\#(y_i)$, replace each occurrence of $c_i$ in $\phi$ by $y_i$
and add the conjuncts $\exists y_i.\alpha_i^\#(y_i)$
(as bounded quantifiers)  to the formula so obtained.
(ii) follows from commutativity of $F$ and then
the first claim in (iii) by disjoint normal form
of Boolean expressions.  Now assume that
 the $\alpha_i$ are in $F_0$ whence  defined within $F$ by
$s_i \alpha_i =t_i$ with constant $\Lambda_F$-terms $s_i,t_i$.
In this case, multiply  any equation  $p(\bar x,\bar c)=0$
with all the constant $\Lambda_F$-terms  $s_j^{k_j}$
 where $k_j$ is the maximum power to which
$c_j$ occurs in the equation. (iv) is obvious.
\end{proof}

For the following compare
\cite[Fact 3.2]{hz}. Call an equation in $\Lambda_L$ 
\emph{special} if it is of the form  $x=0$, $z=x+y$, or $y=x^\perp$.

\begin{fact}\lab{log}
For every quantifier free formula $\phi(\bar x)$ in $\Lambda_L$
there is a conjunction  $\phi'(\bar x,\bar z, \bar y)$
of special equations
 with new variables $\bar z,\bar y$ 
and a boolean combination $\phi''(\bar x,\bar y)$ of equations between variables 
from $\bar x,\bar y$
such  that $\phi(\bar x)$ 
is equivalent within lattices with involution to both
\[\begin{array}{l}
\phi^\exists(\bar x,\bar z,\bar y) \equiv
\exists \bar z \bar y.\, \phi'(\bar x,\bar z,\bar y) \wedge
\phi''(\bar x,\bar y)\\
\phi^\forall(\bar x,\bar z,\bar y) \equiv
\forall \bar z \bar y.\, \phi'(\bar x,\bar z,\bar y) \Rightarrow
\phi''(\bar x,\bar y).
\end{array}\]
 Moreover, for any $\bar u$ in an involutive lattice $\vb$  there are
unique $\bar v,\bar w$ such that $\vb\models \phi'(\bar u,\bar w,\bar v)$.
Also, if $\phi$ is a conjunction (disjunction) 
of equations and negated equations then so is $\phi''$.
\end{fact}

\section{Review of earlier work}\lab{invar}
Recalling some definitions and results from Sections 8--10 of \cite{hz}, we shall relate $\vb$ and $F$, directly,
without using endomorphism rings as intermediate step.
Given a basis $\bar v$  define
the relation $\theta_{\bar v}$ between
$\vb^n$ and $\fn$  by
\[ \bar u \;\theta_{\bar v} \; \bar A \;\;\mbox{ iff } 
\bar u  =(\Sp_{\bar v}(A_1),\ldots ,\Sp_{\bar v}(A_n))\; \mbox{ for } \bar u \in \vb^n,\,
\bar A  \in \fn\]
where $\Sp_{\bar v}(A)$ denotes $\Sp\{ \sum_{i=1}^d a_{ij}v_i\mid j=1,\ldots ,d\}$; that is, the span of the columns of $A$
if $V$ is identified with $F^d$ via the basis $\bar v$. 
This gives rise to a Galois connection and thus to  maps
$ \theta_{L \bar v  F}$ and $\theta_{F \bar v  L}$
mapping  subsets $M$ of $\vb^n$ to  
subsets $K$ of $\fn$ and vice versa
\[ \theta_{L \bar v  F}(M)=\{\bar A  \in \fn\mid  \exists
\bar u \in \vb^n.\; \bar u \theta_{\bar v } \bar A\}\]\[
\theta_{F \bar v  L}(K)=\{\bar u  \in \vb^n\mid  \exists
\bar A \in \fn.\; \bar u \theta_{\bar v } \bar A\}
\]
Our objective
is to match definable subsets
of $\vb^n$ with definable subsets
of $\fn$  via
suitable translations -- preserving quantification type
as much as possible. Observe that $\fd$ 
is not considered  a ($*$-)ring, but 
just a power of the field $F$, formatted in a suitable fashion. 
Definable subsets have certain invariance properties
which we now recall from \cite{hz}, Section 9.

\medskip
Let 
$F^+$ the multiplicative 
subgroup $\{r\mid 0\neq r= r^* \in F \}$ of $F$ and
$\Ot^+(V)$ consist of all \emph{scaled orthogonal maps} $g$: 
for some  $r\in F^+$ and orthogonal linear transformation $h$
\[ g v=h(vr) \;\mbox{ for all } v \in V.\]
Observe  that  the group $Q^+(V)$ is  not changed 
if the form on $V$ is 
scaled by an  element of $F^+$.

 For any  $g \in \Ot^+(V)$, the
map
  $U \mapsto g^L(U):=g(U)$
is an automorphism of $\vb$ (actually, for $d\geq 3$ any automorphism
of $\vb$ is of this form).
We say that  $M \subseteq \vb^n$  is \emph{invariant}
if it is  invariant under the component-wise
action of the  $g^L$, $g \in \Ot^+(V)$.
Clearly,  definable  $M$ are invariant.

Let $\GL(F,d)$ denote the set of of invertible
 matrices in $F^{d\times d}$.
We define 
 $\Ot^+_{\bar  \alpha}(F,d)$,
to  consist of the $T\in\GL(F,d)$
such that   $T^{\dagger_{\bar  \alpha}} =r T^{-1}$ for some $0\neq r \in F^+$,
and consider the  action $A\mapsto TAT^{-1}$ on $F^{d \times d}$.
We call  $K \subseteq \fn$  ${\bar \alpha}$-\emph{invariant} if it is invariant  under the component-wise action of  
 $\Ot^+_{\bar \alpha}(F,d)$;
 \emph{right invariant}
if $(A_1T_1, \ldots ,A_nT_n) \in K$
for all $(A_1, \ldots ,A_n) \in K$ and $T_i \in \GL(F,d)$;
and $\bar \alpha$-\emph{bi-invariant} if both conditions are satisfied, i.e. 
if 
$(TA_1T_1, \ldots ,TA_nT_n) \in K$
for all $(A_1, \ldots ,A_n) \in K$, $T\in\Ot^+_{\bar \alpha}(F,d)$, and  $T_i \in \GL(F,d)$.

Of course, given a  first order formula, 
 right invariance, of the subset of $(F^{d\times d})^n$
it defines, can be stated by a first order sentence;
similarly $\bar \alpha$-invariance if $\bar \alpha$ is definable.
Observe that 
two matrices span the same space of columns if and only if
one is obtained from the other applying an invertible matrix
on the right.
This yields the following, cf.\cite[Fact 9.2]{hz}.

\begin{fact} \lab{linv} 
 $\theta_{L \bar v F}$ and  $\theta_{F \bar v L}$
induce mutually inverse bijections between the set of all
 subsets $M$ of  $\vb^n$
and the set of all right invariant  subsets $K$ of $\fn$.
\end{fact}

Given $\bar \alpha$, define
$\theta_{L\bar \alpha F}(M)$, and
$\theta_{F\bar \alpha L}(K)$, respectively,
as the union of the 
$\theta_{L\bar v F}(M)$, and
$\theta_{F\bar v L}(K)$
where $\bar v$ ranges over all  $\perp$-bases $\bar v$ with
$|\bar v|=\bar \alpha$.
The following is \cite[Proposition 9.3]{hz}.

\begin{proposition}\lab{f3} 
 Fix a  $\perp$-basis
$\bar v$ and $\bar \alpha=|\bar v|$.
Then $\theta_{L\bar v F}$ and $\theta_{F \bar v L}$
induce mutually inverse bijections between the 
set of all invariant
 $M\subseteq \vb^n$  and the set of all $\bar \alpha$-bi-invariant 
 $K\subseteq \fn$;  moreover, for such $M$ and $K$,
$\theta_{L\bar vF}(M)  = \theta_{L\bar \alpha F}(M)$ and
 $\theta_{F \bar v L}(K)=\theta_{F \bar \alpha L}(K)$.
\end{proposition}

In Sections 8 and 10 (cf. Theorem 10.4(ii)) of \cite{hz} we have constructed
translations $\tau^\exists_{L \bar \alpha F}:\Lambda_L \to \Lambda_F^+$
and $\tau^\exists_{F \bar \alpha L}:\Lambda_F^+ \to \Lambda_L$
which preserve quantification type $\Sigma_1$ 
and such that, for any admissible $\bar \alpha$,
$\tau^\exists_{L \bar \alpha F}(\phi(\bar x))$ 
defines $\theta_{L \bar \alpha F}(M)$ if $\phi(\bar x)$ defines $M$ in $\vb^n$
and, if in addition $\bar \alpha$ in $F_0$, 
then $\tau^\exists_{F \bar \alpha L}(\psi(\bar X))$ 
defines $\theta_{F \bar \alpha L}(K)$ if $\psi(\bar X)$ defines $K$ in $\fn$.

\xc{
\section{Examples of translations}
In the subsequent section
we will describe a  procedure
how to  translate from  configurations in lattices to 
expressions in coordinates and prove  correctness of this procedure.
Here, we give an informal account and illustration.
In view of Proposition \ref{f3} we may assume $V=F^d_{\bar \alpha}$ and
$\bar v$ the canonical basis. Elements $u$ of $\vb(V)$
 are given as $u=\Sp(A)$ with $d\times d$-matrices $A$.
Given a formula $\phi(\bar x)$ on the lattice side,
we have to make sure that
the translation  $\tau(\phi)(\bar X)$
(where $\bar X$ is a tuplet of matrices of variables) is satisfied exactly
by those tuplets $\bar A$ of matrices
which have $\bar u$ satisfying $\phi(\bar x)$ where
$u_i =\Sp(A_i)$. Of course, it suffices to deal with atomic
$\phi(\bar x)$, that is equations $t_1(\bar x)=t_2(\bar x)$
where $t_i(\bar x)$ are terms in the language of lattices with
involution.

The basic idea is to describe, for each term $t(\bar x)$
on the lattice side, an algorithm,
depending only on dimension $d$,
how to compute 
(by Gaussian elimination without division) for each 
input tuple $\bar A$ of matrices 
an   output matrix $B$ in wNF such that
  $\Sp(B)=t(\bar u)$ where $u_i=\Sp(A_i)$.
The form elements $\bar \alpha$ are treated as parameters.
The description of the algorithm is via 
an input tuple $\bar X$ of matrices  of variables and, for each branch of the
computation, a 
basic formula
$\pi(\bar X)$, capturing the underlying sequence of cases,
and an  output matrix $P(\bar X)$ of $\Lambda_F^+$-terms,  capturing
the result of computation, and function $f$
giving  pivot positions of the output wNF.

Now, an equation $t_1(\bar x)=t_2(\bar x)$
translates into a formula 
stating that, for each pair of
items  $(\pi_1,f_1,P_1)$ from $t_1$
and  $(\pi_2,f_2,P_2)$
from $t_2$,
one has 
$f_1=f_2$ and 
 $P_1$ equal $P_2$, up to crosswise
multiplications with the pivot terms.

Of course, the description of calculations will involve a number of cases
 rapidly growing
with $d$. 
 For example, consider $t_1=u_1+u_2$ and $t_2=u_3^\perp$,
matrices $X_k=(x_{kij})_{ij}$ of variables 
and possible  input matrices $A_k=(a_{kij})_{ij}$, $k,i,j\in\{1,2,3\}$.
Evaluating $t_1$ requires to transform
the compound matrix  $(A_1|A_2)$ into wNF. Concerning the possible
calculations and associated distinction of cases, we 
will discuss only a few examples. First, assume that
$a_{kij}=0$ for $k=1,2$ and $j=2,3$. This amounts to the side
condition
\[\bigwedge_{i=1}^3 x_{1i2}=0 \;\wedge\;
\bigwedge_{i=1}^3 x_{1i3}=0 \;\wedge\;
\bigwedge_{i=1}^3 x_{2i2}=0 \;\wedge\;
\bigwedge_{i=1}^3 x_{2i3}=0.\]
 Thus, by column permutation
and omitting zero columns in positions $j>3$ we arrive 
at the following matrix $B_1$ and associated matrix $P_1$ of terms
\[ B_1= \left( \begin{array}{ccc} a_{111}&a_{211}&0\\
a_{121}&a_{221}&0\\
a_{131}&a_{231}&0\\ 
\end{array}\right), \quad
P_1= \left( \begin{array}{ccc} x_{111}&x_{211}&0\\
x_{121}&x_{221}&0\\
x_{131}&x_{231}&0\\ 
\end{array}\right).\]
Computing the wNF for $B_1$, the first distinction of cases 
is whether $a_{211}=0$ or $a_{211}\neq 0$.
We consider the latter case, only, that is, the condition
\[ \pi_2\equiv \pi_1 \wedge x_{211} \neq 0. \] 
Continuing the computation of wNF under this assumption we
obtain $B_2$ arising from
\[P_2=
\left( \begin{array}{ccc} x_{111}&0&0\\
x_{121}&x_{111}x_{221}-x_{211}x_{121}&0\\
x_{131}&x_{111}x_{231}-x_{211}x_{131} &0\\ 
\end{array}\right).\]The next distinction of cases is whether 
$a_{111}a_{221}-a_{211}a_{121}=0$ or not; 
these correspond to
\[ \pi_3\equiv \pi_2 \wedge\; x_{111}x_{221}-x_{211}x_{121}=0
\mbox{ and }
\pi_4\equiv \pi_2\wedge\; x_{111}x_{221}-x_{211}x_{121}\neq 0. \]
Continuing under condition $\pi_3$ 
one has $ a_{111}a_{231}-a_{211}a_{131}=0$ or not;
these correspond to 
\[\pi_5\equiv \pi_3\wedge\; x_{111}x_{231}-x_{211}x_{131}=0
\mbox{ and } \pi_6\equiv \pi_3\wedge\; x_{111}x_{231}-x_{211}x_{131}\neq 0\]
and matrices  transformed into  wNF with  pivot positions
$f_5=(1)$  and $f_6=(1,3)$; moreover, the matrix of terms associated with $\pi_5$ is
\[P_5=
\left( \begin{array}{ccc} x_{111}&0&0\\
x_{121}&0&0\\
x_{131}&0 &0\\ 
\end{array}\right).\]
Under condition $\pi_4$ one obtains  wNF with pivot positions
$f_4=(1,2)$, namely  $B_4$  arising from $P_4$ given as
\[
\left( \begin{array}{ccc}(x_{111}x_{221}-x_{211}x_{121}) x_{111}&0&0\\
0&x_{111}(x_{111}x_{221}-x_{211}x_{121})&0\\
(x_{111}x_{221}-x_{211}x_{121})(x_{131}-x_{121})&x_{111}(x_{111}x_{231}-x_{211}x_{131}) &0\\ 
\end{array}\right).\]
To consider a second example for computing $t_1$  assume
$a_{111}=0$,
$a_{121}\neq 0$, $a_{211}\neq 0$, and $a_{kij}=0$ for $k=1,2$, $j=2,3$.
This amounts to $f_7=(1,2)$,
\[ \pi_7\equiv \pi_1\wedge\; x_{111}=0\; \wedge\; x_{121}\neq 0\wedge \;x_{211} \neq 0,  \]
\[P_7= 
\left( \begin{array}{ccc} x_{121}x_{211}&0&0\\
x_{121}x_{221} &x_{211}x_{121}&0 \\
x_{121}x_{231}&x_{211}x_{131}&0
\end{array}\right).\]
Evaluating  $t_2$, the underlying form matters, so assume
$\alpha=(1,1,1)$.
First, consider  $A_3$ such that
$a_{331}\neq 0$ and $a_{3ij}=0$, else; that is 
\[\pi_8\equiv\; x_{331}\neq 0\; \wedge \bigwedge_{(i,j)\neq (3,1) } a_{3ij}=0. \] 
Fact \ref{nf}(iii) gives $\Sp(A_3)^\perp= \Sp(a_{331}I- A_3)$, $f_8=(1,2)$ and
\[ P_8= 
\left( \begin{array}{ccc} x_{331}&0&0\\
0&x_{331}&0\\
0&0&0
\end{array}\right).\]
Next consider $A_3$ satisfying
\[\pi_9 \equiv\; x_{311}\neq 0 \;\wedge \;x_{332}\neq 0\;
\wedge x_{331}=0\;\wedge\;x_{333}=0\;\wedge\bigwedge_{i=1,2;j=2,3} x_{3ij  }=0.\]  
A wNF of $A_3$ is 
\[ 
\left( \begin{array}{ccc} a_{332}a_{311}&0&0\\
a_{332}a_{321}&0&0\\
0&a_{311}a_{332}&0
\end{array}\right)\]
whence
\[ \Sp(A_3)^\perp= \Sp\left( \begin{array}{ccc} 0&0&-a_{332}a_{321}\\
0&0&a_{311}a_{332}\\
0&0&0
\end{array}\right).\]
In case $a_{321}\neq 0$ one has $\pi_{10}=\pi_9\wedge x_{321}\neq 0$,
 $f_{10}=(1)$, and 
\[ P_{10}=
\left( \begin{array}{ccc} -x_{332}x_{321}&0&0\\
x_{311}x_{332}&0&0\\
0&0&0
\end{array}\right).  \] 
In case $a_{321}= 0$ one has $\pi_{11}=\pi_9\wedge x_{321}= 0$,
 $f_{11}=(2)$, and 
\[ P_{11}=
\left( \begin{array}{ccc} 0&0&0\\
x_{311}x_{332}&0&0\\
0&0&0
\end{array}\right).  \] 
As explained above, the set of matrix triples $(A_1,A_2,A_3)$
satisfying the identity $\Sp(A_1)+\Sp(A_2)=\Sp(A_3)^\perp$ can  be defined
by a conjunction of formulas, each referring
to a pair of formulas, one describing a calculation for $t_1$, the
other one for $t_2$. The formulas 
which relate to our examples are $\pi_4,\pi_5,\pi_6,\pi_7$
on one side, $\pi_8,\pi_{10},\pi_{11}$ on the other.
First we deal with the cases where one gets the same pivot positions for
$t_1$ and $t_2$. Here, the corresponding conjuncts of the 
defining formula are
\[\begin{array}{rcrcl}   
\pi_4\wedge \pi_8 &\Rightarrow& x_{331}P_4&=&
(x_{111}x_{221}-x_{211}x_{121}) x_{111}P_8\\
\pi_7\wedge p_8 &\Rightarrow& x_{331}P_7&=&  x_{121}x_{211}P_8\\
\pi_5 \wedge \pi_{10} &\Rightarrow&  -x_{332}x_{321}P_5&=& x_{111}P_{10}
\end{array} 
\]
For all other combinations we have to state a contradiction;
that is, to include a the conjunct 
$\pi_i\wedge \pi_{j} \Rightarrow \bot$
for $(i,j)$ in  $\{4,5,6,7\}\times \{8,10,11\}$
with the exception of $(i,j)=(4,8),(7,8),(5,10)$.
}

\section{Translation via Gauss}\lab{sgau}
In \cite{hz} the translation from $\Lambda_L$
to  $\Lambda_F^+$  was constructed
with the $*$-ring of endomorphism of $V$ as 
an intermediate structure. While this translation
preserves quantification type $\Sigma_1$,
a translation from $\Lambda_L$
to $\Lambda_F^+$ preserving quantifier freeness 
can be constructed
 based on Gaussian elimination.
Recall that $F_0$ denotes the prime subfield of $F$
and that $\bar c$ are constants in $\Lambda_F^+$
to be interpreted by form constants $\bar \alpha$.

\begin{theorem}\lab{p7} 
For any fixed $d$, 
there is a map $\tau_{L\bar c F}:\Lambda_L \to \Lambda_F^+$ 
 such that, for any 
$d$-dimensional
inner product space $V$ over a field $F$  
with  admissible $\bar \alpha$,
$\theta_{L\bar \alpha F}(M) \subseteq \fn$ is defined
by $\tau_{L\bar \alpha F}(\phi)(\bar X)$
if
 $M \subseteq  \vb(V)^n$ 
is defined by
$\phi(\bar x)$ in $\Lambda_L$;
$\tau_{L\bar c F}(\phi)(\bar X)$ is quantifier free if so is
$\phi(\bar x)$.
Moreover,
if $\bar \alpha \in F_0^d$, then
 $\tau_{L\bar \alpha F}(\phi)(\bar X)$ is 
in $\Lambda_F$.
\end{theorem} 
In particular, the translation $\tau_{L \bar c F}$
is uniform for all $V$ with $\dim V=d$; and
$\tau_{L\bar \alpha F}(\phi)$ is obtained
from $\tau_{L \bar c F}(\phi)$ substituting $\bar \alpha$ for
$\bar c$. 
The proof needs some preparations.

\begin{fact}\lab{gau}   Fix $m$ and 
a $d \times m$-matrix $X$ of variables. 
\begin{enumerate} 
\item {\rm Normal form.} For any   $f\in d^\#$ and $m\leq d$,
 there is a quantifier free   formula $\nu^{mf}(X)$
in $\Lambda_F$
such that for any $A \in \fdm$  one has
$F\models \nu^{mf}(A)$ if and only if
$A$ in $f$-wNF.
\item {\rm Computation of normal forms: distinction of cases.}
With a
$d\times m$-matrix $X$ of variables,
 for each $f \in d^\#$, there is a finite set $\Sigma^{mf}$ 
of basic formulas $\sigma(X)$ in
$\Lambda_F$ 
   and for each $\sigma \in \Sigma^{mf}$ there is a $d \times d$-matrix $P^{\sigma f}(X)=(p_{ij}^{\sigma f})_{ij}$
of   $\Lambda_F$-terms in variables from $X$   such that for any $A \in \fdm$ 
  one has
\begin{enumerate}  
\item $F\models \sigma(A)$ for some $f \in d^\#$ and   $\sigma \in
  \Sigma^{mf}$ 
\item
for all $f \in d^\#$ and  $\sigma \in \Sigma^{mf}$, if 
$F\models \sigma(A)$ then  
$P^{\sigma f}(A)$ is in $f$-wNF  and $\Sp(P^{\sigma f}(A)) =\Sp(A)$.   
\end{enumerate}  
\item {\rm Orthogonals.} 
With  a $d\times d$-matrix $X$ of variables,
 for each   $f\in d^\#$ there is a matrix $Q^f(X)=(q_{ij}^f)_{ij}$ of 
$\Lambda_F^+$-terms in variables from $X$
such that for each admissible $\bar \alpha$,
substituted for $\bar y$, and
 each $A \in \fd$ in $f$-wNF one has 
$\Sp(Q^f(A))=\Sp(A)^\perp$ in $F^d_{\bar \alpha}$.
\end{enumerate} 
\end{fact}
\begin{proof}
(i). The construction of $\nu^{mf}$ is obvious,
(ii) Given  $f \in d^\#$, the  $\sigma \in \Sigma^{mf}$ capture
the distinction of cases  in (column)  Gauss-elimination, applied 
to matrices $A$, to  yield $f$-wNF,  and each $\sigma$ grants
that such exists. This is easily (and tediously) expressed via
basic  formulas. The terms in $P^{\sigma f}$
then combine  the elimination calculations,  followed by 
multiplications with  terms obtained for the pivots.   
(iii) is immediate by (iii) of Fact~\ref{nf}. 
\end{proof}

Recall that $\theta_{L\bar \alpha F}=\theta_{L\bar v F}$ 
where $\bar v$ is any $\perp$-basis with $|\bar v|=\bar \alpha$
(Proposition \ref{f3}). 
Thus, to verify that the translation $\phi \to \tau_{L\bar \alpha F}(\phi)$
matches $M$ with $\theta_{L \bar \alpha F}(M)$,
we may argue based on an unspecified
 $\bar v$; this allows to  identify
$V$ with $F^d_{\bar  \alpha}$.
We now have to explain how to relate $\Lambda_L$-terms
 to terms and quantifier free formulas in $\Lambda_F^+$. We associate with each 
variable $x_k$ a matrix $X_k$ of variables.
The following captures the matrix computations 
associated with the evaluation of $\Lambda_L$-terms.

\begin{lemma}\lab{claim}
One can associate with
each $\Lambda_L$-term
 $t(\bar x)$
a finite set $\Gamma_t$ of pairs $(\pi(\bar X),f)$,
where $\pi(\bar X)$ is  a basic formula in 
$\Lambda_F^+$  and $f\in d^{\#}$,
and with each $(\pi(\bar X),f)\in \Gamma_t$ a
$d\times d$-matrix $P^{t \pi f}(\bar X)=(p^{t\pi f}_{ij}(\bar X))_{ij}$
of $\Lambda_F^+$-terms,
 such that the following hold in $F_{\bar \alpha}^d$
for any $F$ and admissible $\bar \alpha$: 
\begin{itemize}
\item[(1)]  For any $\bar A \in \fn$
there is  $(\pi(\bar  X),f)\in \Gamma_t$
such that $F\models \pi(\bar A)$
\item[(2)] For any $(\pi(\bar X),f)\in \Gamma_t$
and  any $\bar A \in \fn$, 
if $F\models \pi(\bar A)$ then  
$P^{t\pi f}(\bar A)$ is in $f$-wNF  and 
 $\Sp( P^{t\pi f}(\bar A))  =t(\bar U)$ where  $U_k=\Sp(A_k)$.
\end{itemize}
\end{lemma}
\begin{proof}
For a variable $x$ we choose  $\Gamma_x=\{(\sigma(X),f)
\mid f \in d^{\#}, \sigma(X) \in \Sigma^{df}\}$ and 
   $P^{x\sigma f}(X)$ as in Fact \ref{gau}(ii). 

Given matrices $A,B$,
one obtains $U=\Sp(A)+\Sp(B)$ as 
$\Sp(C)$ with $C$ derived
as wNF  from the compound matrix
$(A|B)$, omitting the last $d$ zero columns. 
Let $m=2d$ and denote by $(X|Y)$ the $d \times m$-matrix obtained from
the $d \times d$-matrices
$X,Y$ of variables. By
 Fact~\ref{gau}(ii)
one 
 obtains 
  a finite set $\Gamma$ of pairs 
$(\sigma(X|Y),f)$, where $\sigma(X|Y)$ is a
 basic formula in $\Lambda_F$ and $f \in d^{\#}$, 
and for each $(\sigma,f)
 \in \Gamma$  a $d \times d$- matrix $S^{\sigma f}(X|Y)$ of
 terms in $\Lambda_F$
 such that   
for any $A,B \in F^{d \times d}$ there is $(\sigma,f) \in \Gamma$ with 
 $F\models \sigma(A|B)$ 
and 
$S^{\sigma f}(A|B)$ in $f$-wNF
and with
$ \Sp(A)+\Sp(B)= \Sp(S^{\sigma f}(A|B))$. 
Now,
for  $t(\bar x)=t_1(\bar x)+t_2(\bar x)$ let 
 $\Gamma_t$ consist of all pairs $(\pi(\bar X),f)$ where
$(\sigma,f) \in \Gamma$ and
 \[\pi(\bar X)\equiv
\pi_1(\bar X)\wedge \pi_2(\bar X) \wedge \sigma 
(P^{t_1\pi_1 f_1}(\bar X)|P^{t_2\pi_2 f_2}(\bar X))\]
with  $(\pi_i(\bar X),f_i) \in \Gamma_{t_i}$. Put
\[ P^{t\pi f}= S^{\sigma f}(P^{t_1\pi_1 f_1}(\bar X)|P^{t_2\pi_2 f_2}(\bar X)). \]
By Fact~\ref{gau}(iii),
for  any $A \in F^{d \times d}$ in $g$-wNF,
one  obtains $U=\Sp(A)^\perp$ as 
$\Sp(C)$ with $C=Q^g(A)$ and can apply Fact~\ref{gau}(ii)
to transform $C$ into wNF. Formally, this proceeds as follows.
Again, Fact~\ref{gau}(ii)
yields for each $g\in d^{\#}$ a finite set $\Gamma_g$ of pairs $(\sigma,f)$
and for each 
$(\sigma,f)\in \Gamma_g$
a matrix $P^{\sigma fg}(X)$ of $\Lambda_F$-terms 
such that for any $A \in F^{d \times d}$ in 
g-wNF
 there is $(\sigma,f) \in \Gamma_g$ with 
 $F\models \sigma(Q^g(A))$ and  
 in that case  
$R^{\sigma f g}(A)$
is in $f$-wNF  and $\Sp(R^{\sigma f g}(A))=\Sp(A)^\perp$
where $R^{\sigma f g}(X)=P^{\sigma f }(Q^g(X)$.
 Now, for $t(\bar x)=t_1(\bar x)^\perp$ let 
 $\Gamma_t$ consist of all
$(\pi(\bar X),f)$ where $(\sigma,f) \in \Gamma_g$ and
 \[\pi(\bar X)\equiv
\pi_1(\bar X) \wedge \sigma 
(P^{t_1\pi_1 g}(\bar X))\]
with $(\pi_1(\bar X),g) \in \Gamma_{t_1}$. Put
\[ P^{t\pi f}= R^{\sigma fg}(P^{t_1\pi_1 g}(\bar X)). \]
This provides the translation of
 $\Lambda_L$-terms $t(\bar x)$.
\end{proof} 

\begin{proof} of Thm.\ref{p7}.
To deal with equations,
in view of Lemma~\ref{claim},     define $\gamma_{t_1t_2}(\bar X)$ as
the conjunction of all the following  implications
where $(\pi_i,f_i) \in \Gamma_{t_i}$
and  $p_0^{t_i\pi_i f}$ the entry in position  $(f(1),1)$ of 
$ P^{t_i\pi_i f}(\bar X)$: 
 \[ 
 \pi_1(\bar X) \wedge \pi_2(\bar X ) \;\Rightarrow\;
 p_0^{t_2\pi_2 f}(\bar X)P^{t_1\pi_1 f}(\bar X)=
 p_0^{t_1\pi_1 f}(\bar X)P^{t_2\pi_2 f}(\bar X)\] where $f_1=f_2=f$ and
\[  \pi_1(\bar X) \wedge \pi_2(\bar X ) \;\Rightarrow\; x_{111}\neq x_{111}\]
where $f_1 \neq f_2$.
That is, this formula expresses that for 
any substitution $\bar A$ for $\bar X$,
the evaluation of $t_1$ and $t_2$ (according to the
relevant  distinction of cases) yields the same matrix in wNF
up to crosswise multiplying with the first pivots. 
Thus, for
 all $\bar A$  and $U_k=\Sp(A_k)$
\[t_1(\bar U)=t_2(\bar U) \mbox{ if and only if } F\models \gamma_{t_1t_2}(\bar A).\]
The $\gamma_{t_1t_2}(\bar X)$ 
 give the required translations to $\Lambda_F^+$ for 
equations $t_1(\bar x)=t_2(\bar x)$, that is, the
atomic 
$\Lambda_L$-formulas. This  then extends, canonically, to   quantifier free formulas
and further to prenex formulas. Fact~\ref{def} yields the last claim
of the theorem.
 \end{proof}

Theorem \ref{p7} modifies to the case of dimension restrictions.
Given $\ul{d}=(d_1,\ldots,d_n)$, $d_k\leq d$,
define $\vb^{\ul{d}}=\{\bar u \in \vb^d\mid \dim u_k=d_k\}$.
Define $\delta_{h}(X)$ as the
quantifier free formula which is the disjunction of all
$\sigma(X)
\in\Sigma^{df}(X)$ from (ii) in Fact~\ref{gau}
where $|f|=h$.

\begin{corollary}\lab{p7c}
Under the hypotheses of Theorem \ref{p7},
$\theta_{L \bar \alpha F}(M\cap \vb^{\ul{d}})$
is defined by $\tau_{L \bar  \alpha F}(\phi)(\bar X) \wedge \bigwedge_{k=1}^n 
\delta_{d_k}(X_k)$. In particular, if $N$
is quantifier free definable
relatively to  $\vb^{\ul{d}}$
then $\theta_{L \bar \alpha F}(N)$
is quantifier free definable within $(F^{d \times d})^n$.
\end{corollary}

\begin{corollary}
Fix $d$. For any universal sentence $\phi\in\Lambda_L$
there is a universal sentence $\hat{\phi}\in \Lambda_F^+$
such that, for any $d$-dimensional inner product space
$V$ over a field $F$ with  admissible $\bar \alpha$,
$\vb(V)\models \phi$ if and only if 
$F_\alpha\models \hat{\phi }$.
$($Here, $F_\alpha$ has the parameters $\bar c$ interpreted as $\bar \alpha$.$)$
 The analogue holds with
dimension restrictions on variables in $\phi$. 
\end{corollary}

\begin{proof}
Consider $\phi\equiv \forall \bar x.\,\psi(\bar x)$
with quantifier free  $\psi(\bar x)$.
The following $\hat{\phi}$ will do:
$\forall\bar X.\; \tau_{L\bar \alpha F}(\neg\psi)(\bar X)\Rightarrow 0=1$. 
\end{proof}

\section{The Grassmann-Pl\"ucker point of view}\lab{gras}
An alternative to describing subspaces via matrices is
to use  Pl\"ucker coordinates
 (cmp. \cite[Ch. VII]{hodge}, \cite[\S14.1]{Miller}).  
Fix $d$
and consider
 $0<k\leq m\leq d$. 
Let $\mc{ F}_k$ be the set of all pivot positions $f$ 
for  $d \times m$-matrices $A$  in wNF with $\rk(A)=k$,
that is, the set of all strictly monotone maps $M_k \to M_d$.
One has a quantifier free 
formula $\rho^k(X)$ in $\Lambda_F$ such that $F\models \rho^k(A)$ if and only if
$\rk(A)=k$: requiring all  $\ell\times \ell$-subdeterminants  to be zero if
 $\ell >k$ but some to be non-zero for  $\ell =k$.

Let  $I_k$ denote  the set  of all 
 $k$-tuplets $\bar i=(i_1<i_2< \ldots i_k)$,
ordered lexicographically.
 For a $d \times m$ matrix $X$ 
of variables
and $\bar j,\bar i \in  I_K$ 
 let 
 $X^{\bar j}$ the $d \times k$-matrix
where  the $h$-column is the $j_h$-column of $X$ 
and
$X_{\bar i} ^{\bar j}$ 
the $k\times k$-matrix  where the $\ell$-th row
is the $i_\ell$-th row of $X^{\bar j}$. Define the following terms and
$|I_k|$-tuplets of terms in $\Lambda_F$
\[ \pd_{\bar i}^{\bar j}(X) =  \det(X_{\bar i}^{\bar j})\mbox{ and }
 \pd^{\bar j}(X)= (\pd_{\bar i}^{\bar j}(X)\mid \bar i \in I_k);\]
\[ \pd_{\bar i}(X)=\pd_{\bar i}^{\bar j}(X)\mbox{ and }
\pd_k(X)=\pd^{\bar j}(X)\; \mbox{ where } \bar j=(1,\ldots ,k). 
  \]
Then for $\rk(A)=k$ the set
of \emph{projective Grassmann-Pl\"ucker coordinates}
depends on $U=\Sp(A)$, only:  for any $\lambda \neq 0$
\[\pl_k^d(U):= \pl_k^d(A) := \Sp\{ \pd^{\bar j} (A) \mid  
\bar j \in I_k\}\setminus 0 =\pd_k(\lambda A^\#)F \setminus 0.\]

\begin{fact} \lab{uniq}
Let $A,B\in F^{d \times m}$ and $\rk(A)=\rk(B)=k>0$. Then $\pl^d_k(A) \cap \pl^d_k(B)\neq 0$
if and only if $\pl^d_k(A) =\pl^d_k(B)\neq 0$
if and only if  $B=AT$ for some $T \in \GL(F,m)$
\end{fact}
\begin{proof} The first equivalence follows from 
$\dim \pl^d_k(A) = \dim \pl^d_k(B)=1$;
the second  from the fact that $\pl^d_k(A)$
depends only on $\Sp(A)$.
\end{proof}

Recall that, for $0<k\leq d$,
$\Gamma_k^d(F):=\bigcup\{\pl_k^d(A)\mid A \in \fd, \rk(A)=k\}$  
is defined within $F^{I_k} \setminus 0$ 
by the conjunction  
 of the  homogeneous   \emph{Pl\"ucker relations}
in  variables $\ul{y}=(y_{\bar i}\mid \bar i \in I_k)$
and thus may be considered  a projective variety, 
 the \emph{Grassmannian} $\tilde{\Gamma}^d_k(F)$.
 To deal with $k=0$, 
put $\pd_0(X)=0$ (the constant of $\Lambda_F$) if $X$ is the unique
(empty) $d \times 0$-matrix and put  $\pl_0^d(0)=0$ and 
$\Gamma_0^d(F)=0$.

\begin{lemma}\lab{hermite}
{\rm Normal form from Grassmann-Pl\"ucker coordinates.} 
 For each $0<k \leq d$, and $f \in    \mc{ F}_k$
there are a quantifier free  $\Lambda_F$-formula 
  $\pi_f(\ul{y})$
in
variables $(y_{\bar i}\mid \bar i \in I_k)$  
and terms $p_0^f(\ul{y})$ and
a matrix   $P^f(\ul{y})=(p_{ij}^f(\ul{y}))_{ij}$  of terms,
all in $\Lambda_F$, 
 such that for any field $F$ and  $\ul{r} \in\Gamma^d_k(F)$ one has  
$F\models \pi_f(\ul{r})$ if and only if
$p_0^f(\ul{r}) \neq 0$ and  
$\lambda \ul{r} =   \pd_k(A)$  for some $\lambda \neq 0$ and
 matrix $A\in F^{d \times k}$ in  $f$-wNF, 
namely $\lambda= p_0^f(\ul{r})$ and  $A=(p_{ij}^f(\ul{r}))_{ij}$.
\end{lemma} 
\begin{proof}  Cf. the proof of  Theorem II in Chapter VII
of \cite{hodge}. 
For  $f \in \mc{ F}_k$ let  $\pi_f(\ul{y})$ the
formula with states that 
 the first $\bar i$  in $I_k$  such that
$y_{\bar i} \neq 0$ is
 \[f_0:=(f(1), \ldots ,f(k)).\]
Thus, for a matrix $A$ in wNF  and $\ul{r} \in \pl_k^d(A)$ one has
$F\models \pi_f(\ul{r})$ if and only if $A$ has positions $f$ of pivots.

To define the required terms,
let $I_f$ consist of  all $(i,j)$, $i,j \leq d$ such that
     $f(h)<i <f(h+1)$ and $j \leq h$ 
for some $h$   (where  $f(k+1):=d+1$)
and   put 
\[ f_{ij}=(f(1), \ldots ,f(j-1),f(j+1),
\ldots f(h),i, f(h+1), \ldots, f(k)) \mbox{ if  } j < h\]
\[ f_{ij}=(f(1), \ldots ,f(j-1),i, f(h+1), \ldots, f(k))
\mbox{ if } j=h \]
Now, put   $p_{0}^f(\ul{y})=y_{f_0}^{k-1}$
\[ p_{ij}^f(\ul{y}) = \left\{\begin{array}{ll}
 y_{f_0}  & \mbox{ if } i=f(j)    \\ 
(-1)^{h-j}  y_{f_{ij}}   & \mbox{ if } (i,j)  \in I_f \mbox{ and } 
f(h) <i<f(h+1)\\
0& \mbox{ else } 
\end{array} \right. \]
If  $A$ is in NF with positions $f$ of pivots
 then 
 $A$ is recovered from $\ul{r}=\pd_k(A)$  as $P^f(\ul{r})$,
as required, and $\lambda=p_0^f(\ul{r})=1$.

Now,  assume $\ul{r}\in\Gamma^d_k(F)$ and
$F \models \pi_f(\ul{r})$, in particular $r_{f_0} \neq 0$.
Then $\mu\ul{r}=\ul{s}:=
\pd_k(B)$ for some $\mu \neq 0$ and  $B$ 
which may be chosen in NF.
As $\pi_f$ is ``homogeneous'', one has $F\models\pi_f(\ul{s})$
and $B$ with pivot positions $f$.
It follows $s_{f_0}=\pd_{f_0}(B)=1$ and $\mu=r_{f_0}^{-1}$.
 As observed, above, $B=P^f(\ul{s})$.
We are to determine $\nu$ such that $A=\nu B$
 is as  required in the Lemma.
First we should have $A=P^f(\nu\ul{s})$
since all terms $p^f_{ij}(\ul{y})$ are linear in $\ul{y}$.
Second,  $\pd_k(A)=\nu^k\pd_k(B)=\nu^k\ul{s}$
since we deal with $k\times k$-subdeterminants.
Now, the required 
$\lambda$ is  $\lambda=p^f_0(\bar r)=r_{f_0}^{k-1}$;
 with   $\nu=r_{f_0}$ one obtains, indeed,
 $\lambda \ul{r}=\lambda r_{f_0} \ul{s}=
\nu^k\ul{s}$.
\end{proof}

 Given a  dimension vector $\ul{d}$, define
\[
\Gamma_{\ul{d}}^d(F):=\Gamma_{d_1}^{d}(F)\times \ldots \times \Gamma^{d}_{d_n}(F)\]
and let $\fdn$ consist of all $\bar A\in \fn$ with $\rk(A_k)=d_k$
for $k=1,\ldots ,n$. Observe that
$\Gamma^d_{\ul{d}}(F)$ is defined
within $F^{d_1}\times \ldots \times F^{d_n}$
 by the conjunction of the Pl\"ucker relations
together with $\ul{y}_k\neq 0$ for $d_k\neq 0$ while
 $\ul{y}_k$ is the constant $0$ for $d_k=0$.
On the other hand, $\fdn$
is defined within $\fn$ by the (quantifier free) conjunction of the
$\rho^{d_k}(X_k)$.
Call $\Delta \subseteq \Gamma_{\ul{d}}^d(F)$
  \emph{scalar invariant} 
if \[(\ul{r}_1, \ldots, \ul{r}_n)\in \Delta
\Rightarrow (\lambda_1\ul{r}_1, \ldots, \lambda_n\ul{r}_n)\in \Delta\]
for all $\lambda_k \neq 0$.
Such $\Delta$ may be considered a subset of the product
of Grassmannians  $\tilde{\Gamma}^d_{d_i}(F)$.

Call $\Delta$ $\bar \alpha$-\emph{bi-invariant} 
if $\Delta$ is scalar invariant  
and if, in addition,  for all $T\in \Ot^+_{\bar  \alpha}(F,d)$,  
 $\bar A \in \fdn$,
$\ul{r}_k \in \pl^d_k(A_k)$, and 
$\ul{s}_k \in \pl^d_k(TA_k)$,
 $k=1,\ldots ,n$, one has
\[ (\ul{r}_1, \ldots ,\ul{r}_n)\in \Delta \Rightarrow
(\ul{s}_1, \ldots ,\ul{s}_n)\in \Delta.\]
 Define
\[\theta_{F\dv \Gamma}:\fdn \to \Gamma^d_{\dv}(F),\,\;
\theta_{F\dv \Gamma}(\bar A)=\pl^d_{d1}(A_1)\times\ldots \times \pl^d_{dn}(A_n)\] 
\[\theta_{\Gamma\dv F}:\Gamma^d_{\dv}\to \fdn,\;\;
\theta_{\Gamma\dv F}((\ul{r}_1, \ldots ,\ul{r}_n))=
\{\bar A \in \fdn\mid \forall k\;\ul{r}_k \in \pl^d_{d_k}(A_k)\}.
\]
Use the same notation for the associated maps
from  sets to sets -- taking unions of the images.

 To deal with 
definability, for $d_k>0$
associate with $d\times d$-matrices $X_k$ of
 $\Lambda_F^+$-variables 
a $d_k$-tuplet
 $\ul{y}_k=(y_{k\bar i}\mid \bar i \in I_k)$
of $\Lambda_F^+$-variables, and vice versa;
in case $d_k=0$ we match the constant zero matrix in $\fd$ 
with the constant $0 \in F$.
Also,  with a $\Lambda_F^+$-formula $\psi(X_1, \ldots ,X_n)$ 
 we associate a $\Lambda_F^+$-formula  
 $\chi(\ul{y}_1, \ldots, \ul{y}_n)$, and  vice versa.
Namely, given  $\psi$ 
 choose $\chi \equiv \tau_{F \dv \Gamma}(\psi)$ as follows,
using the conjunction $\eta^d_{d_k}(\ul{y}_k)$
 of Pl\"ucker relations defining $\Gamma^d_{d_k}(F)$ and the
formulas $\pi_f$ and matrices $P^f$ of terms from Lemma~\ref{hermite},  
\[ \bigwedge_{k=1}^n \eta^d_{d_k}(\ul{y}_k)
\wedge \bigl(\bigvee_{f_1 \in \mc{ F}_{d_1}, \ldots,f_n \in \mc{ F}_{d_n}} 
\bigwedge_{k=1}^n \pi_{f_k}(\ul{y}_k) \wedge  \psi(P^{f_1}(\ul{y}_1), \ldots, P^{f_n}(\ul{y}_n))\bigr) 
.\]
Given $\chi$  choose  $\psi \equiv \tau_{\Gamma \dv F}(\chi)$ as 
\[ 
\bigvee_{\bar j_1 \in I_{d_1}, \ldots ,\bar j_n \in I_{d_n}}
\bigwedge_{k=1}^n 
 \rho^{d_k}(X_k^{\bar j_k})
\wedge \chi(\pd_{d_1}(X_1^{\bar j_1}), \ldots, \pd_{d_n}(X_n^{\bar j_n})) 
.\]
Call $\chi$ \emph{scalar invariant} (w.r.t. $F$)
if for all $\ul{r}_i \in  \Gamma^d_{d_i}(F)$, $i=1,\ldots n$,
and all $\lambda_1,\ldots ,\lambda_n$ in $F\setminus\{0\}$ 
one has \[F\models 
 \chi(\ul{r}_1, \ldots, \ul{r}_n) \mbox{ iff }
F\models 
 \chi(\lambda_1\ul{r}_1, \ldots, \lambda_n\ul{r}_n);
\]
that is, iff $\chi$ defines  scalar invariant $\Delta \subseteq
\Gamma^d_{\ul{d}}(F)$.

\begin{lemma} \lab{ff}
$\theta_{F\dv \Gamma}$ and $\theta_{\Gamma\dv F}$
establish mutually inverse bijections
between the set of all right invariant  
subsets $K$ of $\fdn$ and the set of all 
scalar invariant
subsets $\Delta$ of $\Gamma^d_\dv(F)$.  Under this correspondence,
\begin{enumerate}
\item $K$ is $\bar \alpha$-bi-invariant if and only if $\Delta$ is $\bar \alpha$-bi-invariant.
\item
If  $K$ is $\bar \alpha$-bi-invariant and  defined by $\psi$
(within $F$) then 
$\Delta$  is  defined by the scalar invariant
formula $\tau_{F \dv \Gamma}(\psi)$. 
\item
If  $\Delta$ is scalar-invariant and  defined by $\chi$
 then  
$K$  is  defined by $\tau_{\Gamma \dv  F}(\chi)$. 
\end{enumerate}
Here, definability is within $F$ by $\Lambda_F^+$-formulas
 with constants $\bar c$ interpreted as $\bar \alpha$. 
Both translations are uniform for all $F$, preserve
quantifier freeness  and
$\Sigma_1$ and  leave  $\Lambda_F$ invariant.
\end{lemma}
\begin{proof}
Obviously, any $\theta_{F\dv\Gamma}(K)$ is scalar invariant
and any $\theta_{\Gamma \dv F}$ is right invariant 
(by Fact~\ref{uniq}).  
Also, recall that  $\pl^d_{d_k}(A_k)$ is
either  $1$-dimensional or zero (if $d_k=0)$. 

Assume that  $\Delta$ is scalar invariant. Then  
one has   $\bar A \in K:=\theta_{\Gamma \dv F}(\Delta) $ 
if and only if  $(\ul{r}_1, \ldots, \ul{r}_n) \in \Delta$ 
whenever  $\ul{r}_k \in \pl^d_{d_k}(A_k)$  for all $k$.
This applies to any  $\bar r=(\ul{r}_1, \ldots, \ul{r}_n)\in \theta_{F \dv
  \Gamma}(K)$ with suitable $\bar A \in K$ 
  to yield $\bar r \in \Delta$. 
This proves $ \theta_{\Gamma \dv F}(\theta_{F \dv\Gamma}(\Delta))=\Delta$.

Assume that $K$ is right invariant and  $\Delta=
\theta_{F \dv\Gamma}(K)$.  Let  $\bar B \in \theta_{F
  \dv\Gamma}(\Delta)$, that is, there is 
$(\ul{r}_1, \ldots, \ul{r}_n) \in \Delta$ such that
$\ul{r}_k \in \pl^d_{d_k}(B_k)$ for all $k$, whence for some
$\bar A \in K$  
also  $\ul{r}_k \in \pl^d_{d_k}(A_k)$    for all $k$.
By Fact~\ref{uniq},  there are $T_k \in \GL(F,d)$ 
such that $B_k =A_kT_k$ whence $\bar B \in K$ by
right invariance. Thus, $\theta_{F \dv\Gamma}(\theta_{\Gamma \dv F}(K))=K$.

(i) is now obvious. 
(ii) and (iii)  follow, immediately,  in view of Lemma~\ref{hermite} and by inspection of the formulas.
\end{proof}

We write $\Gamma^d_\dv(F_{\bar \alpha})$ to refer to the underlying
space $F^d_{\bar \alpha}$.
Now, assume $d\geq 3$ and define $\vb^\dv=\{\bar u \in \vb^n\mid \dim u_k=d_k\}$ and
\[\theta_{L\dv \alpha \Gamma}:\vb^\dv\to \mc{P}(\Gamma^d_\dv(F)),\;\;
\theta_{L\dv \alpha \Gamma}= \theta_{F\dv \Gamma}\circ \theta_{L \alpha F}|\vb^\dv\]
\[ \theta_{L\dv \alpha \Gamma}(M)=\bigcup_{\bar u\in M} 
\theta_{L\dv \alpha \Gamma}(\bar u) \mbox{ for } M \subseteq \vb^\dv.\]
Recall  from \cite[Fact 5.1]{hz} that $\vb^\dv$ is positive primitive definable within $\vb^n$.
 $M \subseteq \vb^\dv$  is \emph{defined}
by  $\phi(\bar x) \in\Lambda_L$ \emph{within} $\vb^{\ul{d}}$  if 
 $M=\{\bar u \in \vb^\dv\mid \vb\models \phi(\bar u)\}$.
 Theorem~\ref{p7}, Lemma~\ref{ff}, and \cite[Theorem 10.4(ii)]{hz} yield, immediately,
the following where
\[ \tau_{L \dv \alpha \Gamma} =\tau_{F \dv \Gamma} \circ \tau_{L \alpha F},\;\;
\tau_{\Gamma \dv \alpha  L}=  \tau_{F L} \circ \tau_{\Gamma \dv \alpha F}.\]

\begin{theorem}\lab{proj}
Assume admissible $\bar \alpha$ and consider a dimension vector $\ul{d}$.
$\theta_{L\dv \alpha\Gamma}$ establishes a bijection from the set of 
invariant subsets   $M$ of $\vb^\dv$ 
onto the set of $\bar \alpha$-bi-invariant subsets  $\Delta$ of
$\Gamma^d_\dv(F)$. Moreover
\begin{enumerate} 
\item If   $M$ is  defined by $\phi$ within $\vb^{\ul{d}}$ then
$\theta_{L \dv \alpha\Gamma}(M)$ is $\bar \alpha$-bi-invariant and
 defined within  $\Gamma^d_\dv(F_{\bar \alpha})$  by the scalar invariant formula  $\tau_{L \alpha \Gamma}(\phi)$ --  
which,
 is in  $\Sigma_1$ resp.  quantifier free   if so is $\phi$
and which is  in $\Lambda_F$ if    $\bar \alpha \in F_0^d$
\item If    $\Delta$ is $\bar \alpha$-bi-invariant and  defined by $\chi$ 
within $\Gamma^d_\dv(F_{\bar \alpha})$
then
$\theta_{\Gamma \dv \alpha L}(\Delta)$ is 
 defined by $\tau_{\Gamma \dv \alpha L}(\chi)$ within $\vb^{\ul{d}}$ --
which is $\Sigma_1$  if  so is $\chi$.
\end{enumerate} 
\end{theorem}

\begin{corollary} Given any $*$-field $F$,
  $\alpha_i=\alpha_i^*\neq 0$ in $F$,
quantifier free definable $M \subseteq \vb(F^d_\alpha)$  and  $d_1, \ldots ,d_n \leq d$,  the set consisting  of
all  Pl\"ucker coordinates
 $(\pl_{d_1}(u_1), \ldots ,\pl_{d_n}(u_n))$, where  $\bar u \in M$ and
 $\dim
 u_i=d_i$ for $i \leq n$, is a
subset of  the product of the  Grassmannians $\tilde{\Gamma}^d_{d_k}(F)$
which is a disjoint union of sets defined by basic
scalar invariant $\Lambda_F^+$-formulas.   
\end{corollary}
This follows by Fact~\ref{def}.

\section{Homogeneous formulas}\lab{shom}
As mentioned, earlier, Pl\"ucker coordinates are projective coordinates;
thus, one should look for ``homogeneous'' descriptions of definable sets.
First,  variables will be sorted according to dimensions $k\leq d$
and come as strings of pairwise distinct members of the same sort, corresponding to the dimension
$\binom{d}{k}$ of  Grassmann-Pl\"ucker
coordinates in dimension $k$; we write $\bar x \in X_k$.
A  $*$-polynomial
   is \emph{homogeneous}
if, after replacing  $x^*$ by $x$ for each variable $x$,
one obtains for each $\bar x \in X_k$
a  homogeneous polynomial in variables 
$\bar x$, considering the others as constants (in a suitable polynomial ring).
An equation in $\Lambda_F$  is \emph{homogeneous} if it is
of the form $p=0$ with homogeneous $*$-polynomial $p$ with integer coefficients;
a formula in $\Lambda_F$  is \emph{homogeneous}
if each of its atomic subformulas is a homogeneous equation.
'Contradiction' $\bot$ and 'tautology' $\top$ will be considered  homogeneous equations.

 \medskip
For $\ul{r} \in \Gamma_k^d(F)$ 
let $\theta_k(\ul{r})$ denote the unique subspace  of $F^d$
such that $\ul{r} \in \pl_k^d(U)$. 
Also, let $f_{\ul{r}}\in \mc{F}_k$ give the pivot positions
of rank $k$ matrices $A$ in wNF, according to Lemma~\ref{hermite},
such that $\pd_k(A)\in \theta_k(\ul{r})$. From the proof of Lemma~\ref{hermite} we have
\begin{fact}\lab{hom}
 For each $f \in \mc{F}_k$ there are a basic homogeneous formula
$\pi_f(\bar x)$, $\bar x \in X_k$, and a $d\times k$ matrix $A_f(\bar x)$ of integer multiples of the $x_i$
 such that, for any $\ul{r} \in \Gamma_k^d(F)$,
$F\models \pi_f(\ul{r})$ iff $f_{\ul{r}}=f$ and then $\ul{r}=\pd_k(A_f(\ul{r}))$. 
\end{fact}

Considering interpretations of homogeneous formulas we require
that  $\bar x \in X_k$ is mapped onto some $\ul{r} \in \Gamma^d_k(F)$.
Fixing the values $\alpha_i \in F$ for the constants $c_i$,
\emph{validity} of a homogeneous formula $\phi(\bar x_1,\ldots ,\bar x_n)$
under $\bar x_i \mapsto \ul{r}_i$
 has an obvious meaning  in the atomic case, and so in general:
we write $\Gamma^d(F_{\bar \alpha})\models \phi(\ul{r_1}, \ldots ,\ul{r}_n)$. 
To be more formal, $\Gamma^d(F_{\bar \alpha})$
is the multi-sorted structure with sorts $\Gamma^d_k(F)$
and multi-sorted relations given by the homogeneous equations. 

For simplicity, in this section
we assume $V=F^d$ with canonical basis $\bar v$
which is an $\bar \alpha$-basis. Thus, $\vb=\vb(F_{\bar \alpha}^d)$.
For $M\subseteq \vb^n$ and a dimension vector $\ul{d}=(d_1,\ldots ,d_n)$ define 
\[\theta_{\ul{d}}(M)=\{(\ul{r}_1,\ldots ,\ul{r}_n)\mid \ul{r}_i\in \Gamma^d_{d_i}(F),
(\theta_{d_1}(\ul{r}_1),\ldots ,\theta_{d_n}(\ul{r}_n))\in M\}
.\]
Define $\theta(M)\subseteq \Gamma^d(F_{\bar \alpha})$ as the union of the $\theta_{\ul{d}}(M)$ where
$\ul{d}$ ranges over all dimension vectors of length $n$.

\begin{theorem}\lab{homth}
Assume admissible $\bar \alpha \in F_0^d$.
If $M$ is defined within $\vb(F_{\bar \alpha}^d)^{\ul{d}}$
$[$within $\vb(F_{\bar \alpha}^d)^{{n}}$, respectively$]$ by a formula $\phi(\bar x)$ in $\Sigma_k\; (\Pi_k)$, $k \geq 1$,
then $\theta_{\ul{d}}(M)$
$[\theta(M)]$ can be defined within $\Gamma^d_{\ul{d}}(F_{\bar \alpha})$ 
$[$within $\Gamma^d(F_{\bar \alpha})]$ by a homogeneous
$\Lambda_F$-formula in $\Sigma_k\; (\Pi_k)$. 
\end{theorem}

\begin{proof} We prove  special cases, first.
\begin{enumerate} 
\item For any $d_0=d_1\leq d$  there is a conjunction 
$\eta_{d_0d_1}(\bar x_0,\bar x_1)$  of homogeneous equations
(where $\bar x_i \in X_{d_i}$)  such that, for any $\ul{r}_i \in \Gamma_{d_i}^d(F)$,
$\theta_{d_0}(\ul{r}_0)=\theta_{d_1}(\ul{r}_1)$ iff $F\models\eta_{d_0d_1}(\ul{r}_0,\ul{r}_1)$.
\item For any $d_1,d_2 \leq d_0 \leq d$, $d_0\leq d_1+d_2$, 
 there is a conjunction, denoted as
$\sigma_{d_0d_1d_2}(\bar x_0,\bar x_1,\bar x_2)$, 
of homogeneous equations (where $\bar x_i\in X_{d_i}$)
 such that, for any $\ul{r}_i \in \Gamma_{d_i}^{d}(F)$,
one has 
$\theta_{d_0}(\ul{r}_0)=\theta_{d_1}(\ul{r}_1)+\theta_{d_2}(\ul{r}_2)
$ iff $F\models \sigma_{d_0d_1d_2}(\ul{r}_0,\ul{r}_1,\ul{r}_2)$.
\item
Assume admissible $\bar \alpha \in F_0^d$, $d_i\leq d$, $d_0+d_1=d$.
There is a conjunction
$\kappa_{d_0d_1}(\bar x_0,\bar x_1)$ 
of homogeneous equations (where $\bar x_i\in X_{d_i}$)
 such that, for any   $\ul{r}_i \in \Gamma_{d_i}^{d}(F)$,
one has 
$\theta_{d_0}(\ul{r}_0)=(\theta_{d_1}(\ul{r}_1))^\perp
$ iff $F\models \kappa_{d_0d_1}(\ul{r}_0,\ul{r}_1)$.
 \end{enumerate}
In all cases, we  first
consider $f_i \in \mc{F}_{d_i}$
and construct formulas to be applied only to 
$\ul{r}_i \in \Gamma_{d_i}^{d}(F)$ with $f_{\ul{r}_i}=f_i$.
 Form the matrices $A_{f_i}=A_{f_i}(\bar x_i)$
of integer multiples of variables form $\bar x_i$, according to Fact \ref{hom}.
First, assume $d_0<d$.
In (i) and (ii) form the compound  matrices $(A_{f_0}|A_{f_1})$
resp. $(A_{f_0}|A_{f_1}|A_{f_2})$; for  $d_0<d$, 
 the formula $\eta_{f_0f_1}$ in (i) 
states that any $d_0+1 \times d_0+1$-subdeterminant is $0$;
 $\sigma_{f_0f_1f_2}$ in (ii) states that any $d_0+1\times d_0+1$-subdeterminant of  
$(A_{f_0}|A_{f_1}|A_{f_2})$ is $0$ but 
some  $d_0\times d_0$-subdeterminant of  
$(A_{f_1}|A_{f_2})$ is not $0$.
In  (iii)  form the compound matrix $(A_{f_0}|Q^{f_1}(A_{f_1})$
(where $Q^{f_1}(X)$ is as in (iii) of Fact~\ref{gau})
and require by means of $\kappa_{f_0f_1}$   all $d_0+1\times d_0+1$-subdeterminants to be $0$.
Now, put 
\[\sigma_{d_0d_1d_2}\equiv \bigvee_{f_i\in\mc{F}_{d_i},i=0,1,2} \;\bigwedge_{i=0}^2\pi_{f_i}(\bar x_i)
\wedge \sigma_{f_0f_1f_2} (\bar x_0,\bar x_1,\bar x_2)
\]
and similarly for $\eta_{d_0d_1}$ and $\kappa_{d_0d_1}$.
In case $d_0=d$ nothing is required in (i) and (iii), only the 
second part in (ii).

For proving the theorem,
it suffices to consider $\phi$ a conjunction (disjunction)
of equations and negated equations and to
derive a translation in $\Sigma_1$ ($\Pi_1$).
Then $\phi''$ from   Fact~\ref{log} has the same form as $\phi$.
We use $\xi$ and $\xi_i$ as names for variables
occurring  in $\phi'$ or $\phi''$.
Consider  maps $\delta$ associating with 
each $\xi$ a  dimension $\delta(\xi) \in\{0,\ldots ,d\}$.
For each $\delta$ and $\xi$ choose a specific vector 
$\xi^\delta \in X_{\delta(\xi)}$ of variables.
Call $\delta$ \emph{admissible} for $\phi$ if
$\delta(x_i)=d_i$ for all $i$ and if
 for  each
special equation in $\phi'$  the relevant dimension restrictions
are satisfied:
\[\begin{array}{ll}
\delta(\xi_1),\delta(\xi_2) \leq \delta(\xi_0) \leq 
\delta(\xi_1)+\delta(\xi_2)&\mbox{ for } \xi_0=\xi_1+\xi_2\\
\delta(\xi_0)=d- \delta(\xi_1)&\mbox{ for } \xi_0=\xi_1^\perp\\
\delta(\xi_0)=\delta(\xi_1) &\mbox{ for } \xi_0=\xi_1 \\
\delta(\xi_0)=0 &\mbox{ for }  \xi_0=0  
\end{array}
\]
Let $D$ denote the set of all admissible $\delta$.
Observe that any assignment of values $\bar u \in \vb^n$ to $\bar x$
gives rise to  $\delta \in D$: the dimensions
of values $t(\bar u)$ of subterms $t(\bar x)$ under the evaluation of $\phi$.

Given $\delta \in D$, we define the translation $\tau^{Q\delta}$, $Q\in
\{\exists,\forall\}$, first for the special equations
making up  $\phi'$
and let $\tau^{Q\delta}(\phi')$ denote the conjunction of all the \[\begin{array}{lcl}
\tau^{Q\delta}(\xi_0=\xi_1+\xi_2) &\equiv&
\sigma_{\delta(\xi_0)\delta(\xi_1)\delta(\xi_2)}(\xi_0^\delta,\xi_1^\delta,\xi_2^\delta) \\ 
\tau^{Q\delta}(\xi_0=\xi_1^\perp) &\equiv&
\kappa_{\delta(\xi_0)\delta(\xi_1)}(\xi_0^\delta,\xi_1^\delta)  \\
\tau^{Q\delta}(\xi_0=0)&\equiv&\top.
\end{array}
\]
Equalities in $\phi''$ are translated as
\[\tau^{Q\delta}(\xi_0=\xi_1) \equiv \left\{\begin{array}{ll} 
\eta_{\delta(\xi_0)\delta(\xi_1)}(\xi_0^\delta,\xi_1^\delta)&\mbox{ if } \delta(\xi_0)=
\delta(\xi_1)\\ \bot &\mbox{ if } \delta(x_0)\neq \delta(\xi_1).\end{array}\right.  \]
For a negated equality $\beta\equiv\neg(\xi_0=\xi_1)$
 occurring in $\phi''$ we define
$\tau^{Q\delta}(\beta)$ as $\neg \eta_{\delta(\xi_0)\delta(\xi_1)}(\xi_0^\delta,\xi_1^\delta)$
if $\delta(\xi_0)=\delta(\xi_1)$; otherwise,
$\tau^{\exists \delta}(\beta)\equiv \top$
and $\tau^{\forall \delta}(\beta)\equiv \bot$.
Now, $\tau^{\exists \delta}(\phi'')$ is the conjunction of all these,
$\tau^{\forall\delta}(\phi'')$ the disjunction (recall the assumption on
$\phi$ and 
$\phi''$).
With $(x_1,x_2, \ldots)^\delta=(x_1^\delta,x_2^\delta,\ldots)$,
and similarly for $\bar y^\delta$ and $\bar z^\delta$, we finally
arrive at the translations 

\[\begin{array}{lclcl}
\tau^\exists(\phi(\bar x))&\equiv& \bigvee_{\delta\in D} \exists \bar z^\delta\exists 
\bar y^\delta 
\bigl(\tau^{\exists \delta}(\phi')(\bar x^\delta,\bar z^\delta,\bar y^\delta)&\wedge& \tau^{\exists\delta}(\phi'')(\bar x^\delta,\bar y^\delta)\bigr)\\
\tau^\forall(\phi(\bar x))&\equiv& \bigwedge_{\delta\in D} \forall \bar z^\delta 
\forall\bar y^\delta 
\bigl(\tau^{\forall\delta}(\phi')(\bar x^\delta,\bar z^\delta,\bar y^\delta)&\Rightarrow&  \tau^{\forall\delta}(\phi'')(\bar x^\delta,\bar y^\delta)\bigr)\end{array}
\]
For the proof just observe that any substitution for 
the $x_i$ in $\Gamma^d(F_{\bar \alpha})$ gives rise to a  
substitution for the $z_j,y_k$  such that the corresponding
assignment of subspaces   
satisfies $\phi'$; that is, with the associated dimensions
given by $\delta \in $D,
 the  assignment in $\Gamma^d(F_{\bar \alpha})$
satisfies $\tau^{Q\delta}(\phi')$ 
 and  is, in the projective setting,
uniquely determined by the values of the $x_i$.  This 
provides a translation $\tau^{Q\ul{d}}(\phi):=\tau^{Q}(\phi)$
in case of fixed dimension vector $\ul{d}$.
To obtain a defining formula for $\theta(M)$ within $\Gamma^d(F_{\bar \alpha})$,
one just has to form
the disjunction of the $\tau^{Q\ul{d}}(\phi)$
with $\ul{d}$ ranging over all $\ul{d}=(d_1,\ldots ,d_n)$, $d_i\leq d$.
\end{proof}

Recall that $\leq$ may be considered a fundamental relation of
$\Lambda_L$ or defined by $x\leq y \Leftrightarrow x+y=y$.
Modifying (i) and (iii), allowing $d_0 \leq d_1$
respectively $d_0+d_1 \leq d$, one obtains the following.

\begin{corollary} Assume admissible $\bar \alpha \in F_0^d$.
If $M$ is defined within $L(F_{\bar \alpha}^d)^n$  by inequalities of the
form $x_i \leq x_j$ and $x_i \leq x_j^\perp$ then
$\theta_{\ul{d}}(M)$ can be defined within $\Gamma^d_{\ul{d}}(F_{\bar \alpha})$
by a quantifier free homogeneous $\Lambda_F$-formula.
\end{corollary}
Given a partially ordered set $P$ with involution and a space $V$,
an example of such $M$ is obtained as the set 
of all representations of $P$ in $V$. 
Of course, the Corollary applies to $M$ defined by a conjunction
of basic equations. An extension to equations 
given by compound terms appears doubtful;
anyway, the approach of Theorem~\ref{p7}
hardly can be modified  to preserve homogeneity.

\section{Translating lattice formulas}\lab{trans}

Here, we  assume $F$ also
endowed with the  operation $r \mapsto r^{-1}$ where $0^{-1}=0$
and the symbol $^{-1}$ included into the language $\Lambda_F$.
We show that on the lattice side 
quantifier free definability amounts to definability by 
equations $t=0$ and $s=1$. 
For this,
we refer to the concept of frame which 
underlies the coordinatization of modular lattices
(requiring $d\geq 3$):
A   
\emph{frame}  of $\vb=\vb(V)$, $\dim V=d$,   is a system   
 $\bar a=(a_{ij}\mid 1 \leq i,j \leq d)$ 
of elements such that 
for pairwise distinct $i,j,k$ (where $a_i=a_{ii}$)
\[ \One=\bigoplus_\ell a_\ell,\;\;a_{ij}=a_{ji},\;\:  a_i+a_j =a_i\oplus a_{ij},\;\;
a_{ik}=(a_i+a_k)\cap(a_{ij}+a_{jk}).\]
Such  are in correspondence with  bases $\bar v$ via
$a_{ii}=v_iF$, $a_{ij}=(v_j-v_i)F$.
Given a  frame, for any $i\neq j$ there is an isomorphism $\omega^{\bar a}_{ij}$
 of the field
 $F$   onto  $R_{ij}(\bar a)=\{u \in \vb\mid u\oplus a_i=
a_i+a_j\}$ endowed  with  operations
 given by 
 lattice terms with constants from $\bar a$.
For a corresponding basis one has $\omega^{\bar a}_{ij}(r)=(v_j-rv_i)F$.

In
\cite[Lemma 10.2]{hz}, frames $\bar a$ associated with  $\perp$-bases $\bar v$ such that
$\bar \alpha =|\bar v|$ have been characterized
by additional relations. 
 Such frames are called $\bar \alpha$-\emph{frames}.
Here, only $\bar \alpha \in F_0^d$ will be considered.
For all $i\neq j$ one has, besides the orthogonality relations
$a_i\cap a_i^\perp=0$ and $a_i \leq a_j^\perp$,
 the relations $r_{kij}(\bar a)=
(a_i+a_j)\cap a_{ij}^\perp$ where  the 
$r_{kij}(\bar z)$ are lattice terms such
that $\omega^{\bar a}_{ij}(\alpha_k)=r_{kij}(\bar a)$
(observe each element of $F_0$ is given by a constant
$\Lambda_F$-term).
For $\bar \alpha=(1,\ldots ,1)$,
$\bar \alpha$-frames are \emph{orthonormal}
and the additional relations  amount to 
$a_j\ominus_{ij}a_{ij}=(a_i+a_j)\cap a_{ij}^\perp$
where $\ominus_{ij}$ is the term describing subtraction in $R_{ij}(\bar a)$ (cf. proof of \cite[Theorem 2.7]{jacm}.

Given an $\bar \alpha$-frame, $\bar \alpha \in F_0^d$,
there is also a  term in $\Lambda_F$ defining an involution on $R_{21}(\bar a)$
such that $\omega^{\bar a}_{21}$ becomes an isomorphism
of $*$-fields.

\medskip
Given $d\geq 3$,  $\phi(\bar x)\in\Lambda_F$, $\bar x=(x_1,\ldots ,x_n)$,
define $M_\phi=M_\phi(F)$ to consist of all $(\bar a,\bar b)$ in $\vb$
such that $\bar a$ is an $\alpha$-frame and $\bar b \in R_{21}^n$
satisfying $\phi(\bar x)$. Observe that $M_\phi$ consists
of tuplets of $1$-dimensionals, only; let $\ul{d}$
the corresponding dimension vector.

\begin{fact}\lab{frame}
Given an $\alpha$-frame in $\vb$, $\bar \alpha \in F_o^d$, and  quantifier free $\phi(\bar x) \in \Lambda_F$ 
there is a term $t(\bar z,\bar x)$ in $\Lambda_L$
such that $M_\phi=\{(\bar a,\bar b)\mid t(\bar a,\bar b)=0\}$.
\end{fact}
\begin{proof} Observe that a conjunction of equations $t_i=0$
can be comprised to a single one $\sum_i t_i=0$; similarly, $t_i=1$
to $\bigcap_i t_i=1$. Also, 
$t=1$ is equivalent to $t^\perp=0$. Thus, it suffices
to give a list of equations $t_i=0$ and $s_j=1$ which defines $M_\phi$.

First, we have to provide such equations
in variables $z_{ij}$
  defining frames (writing $z_i$ for $z_{ii}$).
Put  $y_{ij}=\sum_{k\neq i,j} z_k$ and $y_i=y_{ii}$. 
That $1=\bigoplus a_\ell$ can be expressed by
$\sum_iz_i=1$ and $\bigcap_i y_i=0$. 
Now, an equation  $a\oplus a_i=a_i+a_j$ 
can be captured by $(x +z_i+z_j)\cap y_{ij}=0$, $x\cap z_i=0$, and $x+y_j=1$;
and  $a=b$, where  $a\oplus a_j=b \oplus a_i= a_i+a_j$,
amounts to  $(a+b)\cap a_i=0$. In particular, 
this applies to the last type of equations defining the concept of
frame and the equations concerning $\bar  \alpha$;
orthogonality  is captured by $a_i \cap (a_i^\perp +a_j)=0$, $i \neq j$.

Recall that  $\phi$ is equivalent in $F$ to 
a disjunction of conjunctions of equations $p_i(\bar x)=0$ and
$q_j(\bar x)\neq 0$  for terms in $\Lambda_F$
 and that any disjunction of equations $t_i=0$ in $\Lambda_L$
can be comprised into a single one (cf. \cite[Fact 6]{hz}).
Now, for any term $p(\bar x)$ in $\Lambda_F$ 
 there  
is a term $\hat{p}(\bar x,\bar z)$ in $\Lambda_L$ having value
$\hat{p}(\bar b,\bar a)=p(\bar b)$ in $R_{21}$ for any
$\alpha$-frame $\bar a$ and   $\bar b\in R_{21}(\bar a)$.
By this, $p_i(\bar x)=0$  can be expressed by 
$(\hat{p}_i(\bar x,\bar z )+z_1)\cap z_2=0$ and $q_i(\bar x)\neq 0$
by $\hat{q}(\bar x,\bar z)\cap z_1=0$.
\end{proof} 
In view of Theorem~\ref{homth} one obtains the following,
improving Theorem 10.4(ii) in \cite{hz} for this special case.

\begin{corollary}\lab{homcor} Let  $d\geq 3$,
$\ul{d}=(1, \ldots ,1)$, and $\bar \alpha \in F_0^d$.
Within $\Gamma^d_{\ul{d}}(F_{\bar \alpha})$, a subset $S$ is definable
by a homogeneous $\Sigma_1$-formula from $\Lambda_F$ if and only if
$S=\theta_{\ul{d}}(M)$ for some $M$ definable within
$(\vb(F)^d)^{\ul{d}}$ by some formula $\exists \bar y.\, t(\bar x,\bar y)=0$.
\end{corollary}

\section{Counterexamples}\lab{count}
Now, we shall give examples where
on the analytic side 
there are restrictions on the possible  descriptions.
For simplicity, we assume
$F$ a $*$-subfield of the complex number field with
conjugation and  $V$ endowed with canonical scalar product
w.r.t. some basis, in particular $\bar  \alpha=(1, \ldots, 1)$. 
In view of this, we may omit reference to $\bar \alpha$
and identify $V$ with $F^d$ -turning the canonical basis $\bar v$
into an $\bar\alpha$-basis. Thus, $\vb=\vb(F_{\bar \alpha}^d)$.
Observe that $\theta_{L\ul{d}\Gamma}=\theta_{F\ul{d}\Gamma}\circ\theta_{L\bar v F}$
by Proposition~\ref{f3}. 
Considering the multi-sorted structure $\Gamma^d(F_{\bar \alpha})$ as in
Section~\ref{shom},
for $M\subseteq \vb^n$  we have $\theta(M)$ the union of
the $\theta_{L\ul{d}\Gamma}$ where $\ul{d}$ ranges over all
dimension vectors of length $n$.
\begin{fact}
Let  $\phi \in \Lambda_L$ and
$\psi\in  \Lambda_F$ be quantifier free formulas such that
$\psi$ defines in $\Gamma^d(F_{\bar \alpha})$ the relation $\theta(M)$
where $M$ is defined in $\vb(F_{\bar \alpha}^d)$ by $\phi$.
Then, for any  
$F'$ a $*$-subfield of $F$, $\phi$ and $\psi$ are related in the same
way.
\end{fact}\lab{sub}
\begin{proof}
Observe that $\vb({F'}^d)$ is embedded into 
$\vb(F_{\bar \alpha}^d)$ by tensoring with $F$ and
the same applies to $\Gamma^d(F_{\bar \alpha}')$ and $\Gamma^d(F_{\bar \alpha})$.
W.r.t. the canonical basis $\bar v$ of ${F'}^d$, the
multi-sorted structure with sorts $\vb({F'}^d)$ and   
$\Gamma^d(F_{\bar \alpha}')$ and relation $\theta$ 
becomes a substructure of the analogous one over $F$.
Also, there is an obvious quantifier free
formula relating the models of $\phi$ and $\psi$  via $\theta$
over $F$. And validity of this formula is inherited by the substructure.   
\end{proof}

The examples use the fact that
with a  formula $\phi$ in $\Lambda_F$
and defining formula $\psi$ for $M_\phi$ 
one can associate a formula $\psi'$ in $\Lambda_F$ which is
equivalent to $\phi$ and inherits structural properties from
$\psi$. This is done as follows.
Put
$v_{ii}=v_i$ and
$v_{ij}=v_j-v_i$ for $i\neq j$,
and $\tilde{v}=(v_{ij}F|i,j\leq d)$;
also, given $\bar r\in F^n$ put  
$\tilde r=((v_1-r_kv_2)F|k\leq n)$.
In view of the isomorphism $\omega_{21}$ it  follows that
$F\models \phi(\bar r)$ if and only if $(\tilde{v},\tilde{r})\in \theta(M_\phi)$. 
Considering any formula $\psi(\tilde{z},\tilde{y})$ defining
$\theta(M_\phi)$, the latter is equivalent
to $F\models \psi(\tilde{v},\tilde{r})$.
Choosing the canonical basis $\bar v$, substituting 
$\tilde{v}$ for $\tilde{z}$
in $\psi$ and, simultaneously,  $\tilde{w}$ for $\tilde{y}$ where
 $w_k=(v_1-x_kv_2)F$, one obtains 
a formula $\psi'(\bar x)$ 
equivalent to $\phi(\bar x)$ in $F$.

\begin{example}\lab{ex} Let $F$ a $*$-subfield of $\mathbb{C}$.
$\theta(M_\phi)$ cannot be defined in $\Gamma^d(F_{\bar \alpha})$ by any formula $\psi$
such that
\begin{enumerate} 
\item  $\psi$ is quantifier free and positive;  here $\phi$ is $x\neq 0$.
\item $\psi$ is a conjunction of equations and negated equations;
here $\phi$ is $x_1=0 \Leftrightarrow x_2 \neq 0$.
\item $\psi$ does not not involve involution; 
here  $\phi$ is $x_2=x_1^*$ and
$F=\mathbb{C}$.
\end{enumerate}
Though, for all these  $\phi$,     $M_\phi$
can be defined within $\vb$ by an equation $t=0$.
\end{example}
\begin{proof}
Definability within $\vb$ follows from Fact~\ref{frame}.
For the negative claims
 consider the associated $\psi'(\bar x)$ and derive contradictions.
In (i) and (ii)
we may assume $F=\mathbb{Q}$ and  $\psi'(\bar x)$ be of the same form as $\psi$
with atomic formulas $p_h(\bar x)=0$, $p_h(\bar x) \in \mathbb{Q}[x]$.
Thus, in (i) $\psi'(x)$
 is built from equations $p_h(x)=0$ by
conjunction and disjunction;
 since $\mathbb{Q}\models \neg\psi'(0)$,
$\psi'(x)$ can have only finitely many satisfying assignments in $\mathbb{Q}$,
in contrast to $\phi(x)$.
In (ii) $\psi'(x_1.x_2)$ would be equivalent to
some   $\bigwedge_{h=1}^mp_h(x_1,x_2)=0\; \wedge q(x_1,x_2)\neq 0$. Observe that for any $r_1\in\mathbb{Q}$
there are infinitely many $r_2$ such that $\mathbb{Q}\models \phi(r_1,r_2)$
which implies that the $p_h$ are zero-polynomials. Thus, $\phi(x_1,x_2)$ would be
equivalent in $\mathbb{Q}$  to $q(x_1,x_2)\neq 0$.
Contradiction, since $\neg \phi(x_1,x_2)$ defines a set which is not closed. 

(iii) Assuming $\theta(M_\phi)$ definable over the field $\mathbb{C}$,
any automorphism of $\mathbb{C}$ would leave $\theta(M_\phi)$
invariant. To arrive at a contradiction,
consider an irreducible $p(x) \in \mathbb{Q}[x]$ of odd degree.
There is an automorphism $\omega$ of $\mathbb{C}$
mapping some  zero $a\in \mathbb{R} $ of $p$ to a zero $b\not\in\mathbb{R}$,
that is $\omega a^*\neq (\omega a)^*$. 
Now consider $\bar a$ the frame given by the canonical basis $\bar v$
and $\bar b= ((v_1-av_1)\mathbb{C},(v_1-a^*v_2)\mathbb{C}))$.
Then $(\bar a,\bar b)$ is in $\theta(M_\phi)$ but its image under $\omega$ is not.  
\end{proof}

In view of Lemma~\ref{ff}, observe the following.
 If  $\psi(\bar X)$
defines bi-invariant $K\subseteq F^{\ul{d}\times \ul{d}}$ where
$\ul{d}=(1, \ldots ,1)\in \mathbb{N}^m$
then $\theta_{F\ul{d}\Gamma}(K)$
is defined within $\Gamma^d_{\ul{d}}(F)=(F^d\setminus\{0\})^m$
by (scalar invariant)  $\chi\equiv\psi(P(\ul{y}_1),\ldots ,P(\ul{y}_m))$
where the $d\times d$-matrix $P(\ul{y})$ has first column $\ul{y}$,
zero else. In particular, if $\psi$ defines $\theta_{L\bar v F}(M)$,
where $M$ consists of $\bar u$ with all $\dim u_i=1$,
 then $\chi$ defines $\theta_{L\ul{d}\Gamma}(M)$
within $\Gamma^d_{\ul{d}}(F)=(F^d\setminus\{0\})^m$,
that is, within the $d-1$-dimensional projective space
over $F$. Thus, the above counterexamples apply
to $\theta_{L\bar v F}$ as well, $\bar v$ the canonical basis.

\section{Absence of involution}
We speak of \emph{absence of involution} if $F$ is just a field,
$V$ a vector space, $\bar v$ any basis,
and
 $\vb$ the lattice 
 of all linear subspaces of $V$, $\Lambda'_L$ and $\Lambda'_F$
the languages of bounded lattices and rings, respectively,
having  multivariate  polynomials   in place of
$*$-polynomials (with integer coefficients).

We consider $\Lambda'_L$ a subset of $\Lambda_L$, generated from variables
by the operations $+$, $0$,  
$t\cap s:=(t^\perp +s^\perp)^\perp$, and  $1:=0^\perp$. 
Observe that for a lattice $\vb$ admitting some
involution, terms from $\Lambda'_L$ can be
evaluated and the value does 
not depend on the choice of involution.
This allows 
 to transfer  results to this case, 
introducing an involution on $\vb$
endowing  $F$ with the identity involution  and  $V$ with a form
declaring some basis orthonormal: that is $\bar  \alpha=(1,\ldots ,1)$.
In this setting, one may read $\Lambda_F$ formulas as such in $\Lambda'_F$,
just omitting $^*$. 

\begin{corollary}
In absence of involution, the following remain valid
(mutatis mutandis):
 Theorem~\ref{p7}, Theorem~\ref{proj}, Theorem~\ref{homth}, and
Corollary~\ref{homcor};
moreover  Fact~\ref{frame} and  Example~\ref{ex} (i), (ii)
 with a conjunction
of equations of the form $t_i=0$, $s_j=1$.
\end{corollary}
 In particular, for $M$ defined by 
 a formula $\phi(\bar x) \in \Lambda'_L$ the 
corresponding $\theta_{L\bar v F}(M)$ is independent of the basis
$\bar v$ and the form $\Phi$ given by declaring $\bar v$  orthonormal.
Indeed, given a second basis $\bar v'$ and the form $\Phi'$ 
rendering $\bar v'$ orthonormal, the linear autmorphism of
the vector space $V$
matching $\bar v$ with $\bar v'$  is also an isometry $\omega$
between these inner product spaces and $\theta_{\bar v}= \theta_{\bar v'}\circ 
\omega$. Finally, observe that  $\omega(M)=M$
since $\phi\in \Lambda'_L$.

Alternatively,  in order to deal with meets in Theorem~\ref{p7},
one may  refer to  the  (column) Zassenhaus algorithm,
 that is, given $A \in F^{d\times k}$ and $B \in F^{d\times \ell}$  
transforming the $2d \times (k+\ell)$-matrix
\[\left(\begin{array}{c|c}A&B\\ \hline A&0
\end{array} \right)
\mbox{ into wNF } \left(\begin{array}{c|c|c}A'&0&0\\ \hline A''&B'&0
\end{array} \right)
 \]
with $A'$ a $d\times m$-matrix, $m=\rk(A)$, and 
$\Sp(B')= \Sp(A) \cap \Sp(B)$.

This gives translations as in Section~\ref{sgau} valid for any choice of basis.
In particular, the formulas and terms
in Lemma~\ref{claim} are in $\Lambda_F$.
In the proof for meets,
we proceed  as for joins, 
with  $d \times d$-matrix  
$M^{\sigma f}(X|Y)$
in $\Lambda_F$ yielding
$\Sp(A) \cap \Sp(B)= \Sp(M^{\sigma f}(A|B))$
 under appropriate  distinction of cases; 
here $M^{\sigma f}$
 describes the calculations in the Zassenhaus algorithm.
To capture the top element of the lattice, one just
refers to matrices having wNF a multiple of $I$.

\medskip
Given a partially ordered set $P$, by Corollary~\ref{homcor}
the set of all representations of $P$ within 
a given vector space $V$ (that is, order preserving maps
$P\to  \vb(V)$) gives rise to a subset of $\Gamma^d(F)$
defined by homogeneous quantifier free formulas.
It remains to clarify how this is related 
to the quiver Grassmannians  (cf. \cite{ring})
of quivers derived from partially ordered sets.

\section{Hidden coordinates}\lab{hid}
Quantifier free translations from $F$ to $\vb$ can be obtained for
quantifier free $\Lambda_F$-formulas which implicitly
provide coordinate systems -- similarly to 
the version  of Cayley factorization due to \cite{Whiteley}.
Assume that $V$ admits an ON-basis,  that is $\bar \alpha=|\bar v|=(1,\ldots,1)$.

We say that $\psi(\bar X)$ in $\Lambda_F$ \emph{hides coordinate systems}
if there are finitely many tuplets
$\bar a^k(\bar x)$ and  $\bar t^k(\bar z)$
of  $\Lambda_L$-terms 
such that, for any $\bar u \in \vb^n$, if $\vb \models
\tau^\exists_{F \bar\alpha L}(\psi)(\bar u)$ then 
there is $k$ such that $\bar a^k=\bar a^k(\bar u)$ 
is an orthonormal frame
 and 
$\bar u=\bar t^k(\bar a^k)$.  

\begin{proposition} Assume $d \geq 3$ and that $V$ admits ON-bases.
If quantifier free  $\psi(\bar X)$ hides coordinate systems 
and defines $\bar \alpha$-bi-invariant $K \subseteq \fn$
then $M=\theta_{F\bar \alpha L}(K)$ can be defined by quantifier free
$\phi(\bar x)$.
\end{proposition}
\begin{proof} 
We  put $\phi(\bar x)\equiv\bigvee_k \phi_k(\bar x)$ where,  
 for any $k$,
$\phi_k(\bar x)$ is the  quantifier free
formula
\[ \bigwedge_{i=1}^n x_i=  t^k(\bar a^k(\bar x)) \;\wedge \tau_{F\bar \alpha L}(\psi)(\bar x,\bar a^k(\bar x)).\]
By hypothesis, if $\vb \models \tau^\exists_{F \bar \alpha L}(\psi)(\bar u)$
then $\vb \models\phi(\bar u)$.
Conversely, if $\vb \models\phi_k(\bar u)$ then 
$\bar a^k(\bar u)$ witnesses
$\exists \bar z$ in  $\vb \models  \tau^\exists_{F \bar \alpha L}(\psi)(\bar u)$.
\end{proof}

In  absence of involution (where 'frames' correspond to bases), 
 one may
find sufficient conditions
for hidden coordinate systems
using  $m$-tuplets $\bar w$ in  $\vb$
which are \emph{associated to frames} via tuplets 
$\bar a(\bar y)$ and $\bar s(\bar z)$ of  lattice terms,
that is,   
   $\bar a(\bar w)$
is a frame of $\vb$ and 
$\bar w= \bar s(\bar a(\bar w))$;
in particular, the sublattice generated by $\bar w$
  is isomorphic to $\vb(F_0^d)$, $F_0$ the prime subfield of $F$.
For $d=3$ and $m=4$, such quadruples  are given by
$4$ points no $3$ of which are collinear,
or by $3$ non-collinear points and a line
incident with none of these, or the duals of such.
More generally,
for fixed $d\geq 3$ and  $m\geq 4$  there are only finitely many 
isomorphism types   of $m$-tuplets  which are
 associated to frames and a finite collection of terms
providing witnesses for  these associations, uniformly for all $F$ (cf. \cite{gelfand,rahm}).  
 Using   these terms, for any fixed $d\geq 3$, one obtains finitely many
tuplets 
$\bar w^\ell(\bar x)$ and  $\bar v^\ell(\bar z)$ of lattice terms
such that
 $\psi(\bar X) \in \Lambda_F$   hides coordinate systems
provided that,   for any $\bar u \in \vb^n$, if $\vb \models
\tau_{F \bar \alpha L}(\psi)(\bar u)$ then 
there is $\ell$ such that $\bar w=\bar w^\ell(\bar u)$ 
is associated to a frame   and
$\bar u=\bar v^\ell(\bar w)$.  

R

\section*{Acknowledgments}
\addcontentsline{toc}{section}{Acknowledgment}
The second author acknowledges co-funding to EU H2020 MSCA IRSES project 731143
by the International Research \& Development Program of 
the Korean Ministry of Science and ICT, grant NRF-2016K1A3A7A03950702.


\end{document}